\documentclass[reqno,12pt]{amsart}%
\usepackage{amsfonts}
\usepackage{amsmath}
\usepackage{amssymb}
\usepackage{graphicx}%
\newtheorem{theorem}{Theorem}
\theoremstyle{plain}

\newtheorem{definition}[theorem]{Definition}
\newtheorem{proposition}[theorem]{Proposition}
\newtheorem{example}[theorem]{Example}
\newtheorem{lemma}[theorem]{Lemma}
\newtheorem{remark}[theorem]{Remark}
\numberwithin{equation}{section}
\renewcommand\endproof{\hfill $\Box$\vskip 0.15in}

\begin{document}
\title[Stochastic Equations with Spatial Noise]{Stochastic Partial Differential Equations Driven by Purely Spatial Noise}
\author{S. V. Lototsky}
\curraddr[S. V. Lototsky]{Department of Mathematics, USC\\
Los Angeles, CA 90089}
\email[S. V. Lototsky]{lototsky@math.usc.edu}
\urladdr{http://math.usc.edu/$\sim$lototsky}
\author{B. L. Rozovskii}
\curraddr[B. L. Rozovskii]{Division of Applied Mathematics \\
Brown University\\
Providence, RI 02912}
\email[B. L. Rozovskii]{rozovsky@dam.brown.edu}
\urladdr{http://www.dam.brown.edu/people/rozovsky.html}
\thanks{S. V. Lototsky acknowledges support from the Sloan Research Fellowship and the
NSF CAREER award DMS-0237724.}
\thanks{B. L. Rozovskii acknowledges support from NSF Grant DMS 0604863, ARO Grant
W911NF-07-1-0044, and ONR Grant N00014-07-1-0044.}
\thanks{An earlier version of the paper was posted on Archives
http://arxiv.org/abs/math.PR/0505551, May 25, 2005.}
\subjclass[2000]{Primary 60H15; Secondary 35R60, 60H40}
\keywords{Generalized Random Elements, Malliavin Calculus, Skorokhod Integral, Wiener
Chaos, Weighted Spaces}

\begin{abstract}
We study bilinear stochastic parabolic and elliptic PDEs driven by purely
spatial white noise. Even the simplest equations driven by this noise often do
not have a square-integrable solution and must be solved in special weighted
spaces. We demonstrate that the Cameron-Martin version of the Wiener chaos
decomposition is an effective tool to study both stationary and evolution
equations driven by space-only noise. The paper presents results about
solvability of such equations in weighted Wiener chaos spaces and studies the
long-time behavior of the solutions of evolution equations with space-only
noise.

\end{abstract}
\maketitle

\today



\section{Introduction}

Stochastic PDEs of the form
\begin{equation}
\dot{u}(t,x)={\mathbf{A}}u(t,x)+{\mathbf{M}}u\left(  t,x\right)  \cdot\dot
{W}\left(  t,x\right)  , \label{I1}%
\end{equation}
where $\mathbf{A}$ and $\mathbf{M}$ are linear partial differential
operators and \ $\dot{W}\left(  t,x\right)  $ is  space-time noise,
are usually referred to as bilinear evolution SPDEs.\footnote{Bilinear SPDEs
differ from linear by the term including multiplicative noise. Bilinear SPDEs
are technically more difficult then linear. On the other hand, multiplicative
models preserve many features of the unperturbed equation, such
as  positivity  of the solution and conservation of mass, and are often more
``physical''.} These
equations are of interest in various applications: nonlinear filtering
for hidden Markov models \cite{Z}, propagation of magnetic field in random
flow \cite{CarMol} stochastic transport \cite{GV, GK, LRps}),
porous media \cite{CarRoz} and others.
The theory and the applications of bilinear SPDEs have been actively
investigated for a few decades now; see, for example, \cite[etc.]{DZ, KR1, NAppl,
Pardoux, Roz, Walsh}.

In contrast, very little is known about bilinear parabolic and elliptic
equations driven by purely spatial Gaussian white noise $\dot{W}\left(  x\right)  .$
Important examples of these equations include:

1. Heat equation with random potential modeled by spatial white noise:%
\begin{equation}
\dot{u}(t,x)=\Delta u(t,x)+u\left(  t,x\right)  \diamond\dot{W}\left(
x\right)  , \label{eq:ranpotent}%
\end{equation}
where $\diamond$\ denotes the Wick product, which, in this case,
coincides with the Skorokhod  integral
in the sense of Malliavin calculus.
A surprising discovery made in \cite{Hu} was that the spatial regularity
of the solution of equation (\ref{eq:ranpotent}) is better than in the case of
similar equation driven by the space-time white noise.

2. Stochastic Poisson equations in random medium:
\begin{equation}
\nabla\left(  A_{\varepsilon}\left(  x\right)  \diamond{\nabla}u\left(
x\right)  \right)  =f(x),
 \label{eq: stochPoisson}%
\end{equation}
\qquad\qquad\ where $A_{\varepsilon}\left(  x\right)  :=\left(
a(x)+\varepsilon{\dot{W}}\left(  x\right)  \right)  ,$ $a(x)$ is a
deterministic positive-definite matrix, and $\varepsilon$ is a positive
number. (Note that $a(x)\diamond{\nabla}u\left(  x\right)  =a(x){\nabla
}u\left(  x\right)  $).

3. Heat equation in random medium:
\begin{equation}
\dot{v}\left(  t,x\right)  =\nabla\left(  A_{\varepsilon}\left(  x\right)
\diamond{\nabla}v\left(  t,x\right)  \right)  +g(t,x).
 \label{eq: stochDiff}%
\end{equation}

Note that the matrix $A_{\varepsilon}$ in equations (\ref{eq: stochPoisson})
and (\ref{eq: stochDiff}) is not necessarily  positive definite; only its
expectation $a\left(  x\right)  $ is.

Equations (\ref{eq: stochPoisson}) and (\ref{eq: stochDiff}) are random
perturbation of the deterministic Poisson and Heat equations. An important
feature of these type of perturbation is that the resulting equations are
\textit{unbiased} in that they preserve the mean dynamics:
the functions $u_{0}\left(  x\right)  :=\mathbb{E}u\left(  x\right)  $ and
$v_{0}:=\mathbb{E}v\left(  t,x\right)  $ solve the deterministic Poisson
equation,
\[
\nabla\left(  a(x){\nabla}u_{0}\left(
t,x\right)  \right)  =\mathbb{E}f\left(  x\right)
\]
and the deterministic Heat equation%
\[
\dot{v}_{0}\left(  t,x\right)  =\nabla\left(  a\left(  x\right)  v_0\left(
t,x\right)  \right)  +\mathbb{E}g(t,x),
\]
respectively.

The objective of this paper is to develop a systematic approach to bilinear
SPDEs driven by purely spatial Gaussian  noise. More specifically, we will
investigate bilinear parabolic equations
\begin{equation}
\frac{\partial v(t,x)}{\partial t}=\mathbf{A}v(t,x)+{\mathbf{M}}%
v(t,x)\diamond{\dot{W}}(x)-f(x) \label{I5}%
\end{equation}
and elliptic equations
\begin{equation}
\mathbf{A}u(x)+{\mathbf{M}}u(x)\diamond{\dot{W}}(x)=f(x).\label{I4}%
\end{equation}
for a wide range of operators $\mathbf{A}$ and $\mathbf{M.}$

Purely spatial white noise is an important type of stationary perturbations.
However, except for elliptic equations with additive random
forcing  \cite{DM, MS, NT}, SPDEs driven by spatial noise have not been
 investigated nearly as
extensively as those driven by strictly temporal or space-time noise.

In the case of spatial white noise, there is no natural and convenient
filtration, especially  in the dimension $d>2$. Therefore, it makes sense to
consider anticipative solutions. This rules out It\^{o} calculus and makes it
necessary to rely on more nuanced techniques provided by Malliavin calculus.

In this paper, we deal with broad classes of operators $\mathbf{A}$ and
$\mathbf{M}$ that were investigated previously for causal (nonanticipating)
solutions of$\mathtt{\ }$the equation (\ref{I1}) driven by space-time white noise.

The notion of ellipticity for SPDEs is more restrictive then in deterministic
theory. Traditionally, nonanticipating solutions of equation (\ref{I1}) were
studied under the following assumptions:

\textbf{(i)}\textbf{ }\textit{The operator} $\mathbf{A}-\frac{1}{2}%
\mathbf{MM}^{\star}$ \textit{is \textquotedblleft elliptic" (possibly
degenerate coercive) operator}.

Of course, this assumption does not hold for equation (\ref{eq: stochDiff})
and other equations in which  the operators $\mathbf{A}$ and $\mathbf{M}$ have
the same order. Therefore, it is important to study equations\textit{
}(\ref{I4}) and (\ref{I5}) under weaker assumptions, for example

\textbf{(ii)}\textit{ The operator }$\mathbf{A}$\textit{ is coercive and }$ord\left(
\mathbf{M}\right)  \leq ord\left(  \mathbf{A}\right)  .$

In 1981, it was shown by Krylov and Rozovskii  \cite{KrRoz} that, unless
assumption (\textbf{i}) holds, equation (\ref{I1}) has no solutions in the
space $L_{2}(\Omega;X)$ of square integrable (in probability) solutions in any
reasonable functional space $X$. The same effect holds for bilinear SPDEs
driven by space only white noise.

Numerous attempts to investigate solutions of stochastic PDEs violating the
stochastic ellipticity conditions were made since then. In particular it was
shown in \cite{LR_AP, MR, NR} that if the operator $\mathbf{A}$ is
coercive ("elliptic") and $ord\left(  \mathbf{M}\right)  \ $is
\textit{strictly }less then $ord\left(  \mathbf{A}\right)  $ then there there
exists a unique generalized (Wiener Chaos) nonanticipative solution of
equation (\ref{I1}).
This generalized solution is a formal Wiener Chaos series  $u=\sum_{\left\vert
\alpha\right\vert <\infty}u_{\alpha}\xi_{\alpha}$, where $\left\{  \xi
_{\alpha}\right\}  _{\left\vert \alpha\right\vert <\infty}$ is the
Cameron-Martin orthonormal basis in the space $L_{2}(\Omega)$. Regularity
of this solution is determined by
system of positive weights $\left\{  r_{\alpha}\right\}  _{\left\vert
\alpha\right\vert <\infty}$ and a function space $X$ such that
\begin{equation}
\left\Vert u\right\Vert _{\mathcal{R},X}^{2}:=\sum_{\left\vert \alpha
\right\vert <\infty}r_{\alpha}^{2}\left\Vert u_{\alpha}\right\Vert
_{L_{2}\left(  (0,T);X\right)  }^{2}<\infty.
 \label{eq: stnorms}%
\end{equation}
The stochastic Fourier coefficients $u_{\alpha}$ satisfy a lower-triangular
system of deterministic PDEs. This system, called {\tt propagator},
 is uniquely determined by the
underlying equation (\ref{I1}).

 Stochastic spaces equipped with the norms similar to (\ref{eq: stnorms})
have been known for quite some time; see, e.g. \cite{kondr1, kondr, NR}. For
historical remarks regarding other types of generalized solutions and
applications to SPDEs see the review paper \cite{LR_shir} and the references therein.

The Wiener Chaos is a bona fide  generalization of the classical It\^{o}
solution: if exists, a non-antisipating square integrable It\^{o} solution
coincides with the Wiener Chaos solution.

In this paper, we establish existence and uniqueness of Wiener Chaos solutions
for stationary (elliptic) equations of the type (\ref{I4}) and evolution
(parabolic) equations of the type (\ref{I5}).
These results are proved under assumption (\textbf{ii}) that allows us in
particular to deal with equations like \eqref{eq: stochPoisson} \ and
\eqref{eq: stochDiff}.
In many cases we are able to find optimal or near-optimal
systems of weights $\left\{  r_{\alpha}\right\}  _{\left\vert \alpha
\right\vert <\infty}$ that guarantee \eqref{eq: stnorms}

Finally, we establish the convergence, as $t\rightarrow+\infty$, of the solution
of the evolution equation to the solution of the related stationary equation.

The structure of the paper is as follows. Section \ref{sec2} reviews the
definition of the Skorokhod integral in the framework of the Malliavin
calculus and shows how to compute the integral using Wiener chaos. Sections
\ref{sec3} and \ref{sec4} deal with existence and uniqueness of solutions to
abstract evolution and stationary equations, respectively, driven by a general
(not necessarily white) spatial Gaussian  noise;
Section \ref{sec4} describes also the
limiting behavior of the solution of the evolution equation; Section
\ref{sec5} illustrates the general results for bilinear SPDEs
driven by purely spatial white noise.

\section{Weighted Wiener Chaos and Malliavin Calculus}

\label{sec2} \setcounter{equation}{0} \setcounter{theorem}{0}

Let ${\mathbb{F}}=(\Omega,{\mathcal{F}},\mathbb{P})$ be a complete probability
space, and ${\mathcal{U}}$, a real separable Hilbert space with inner product
$(\cdot,\cdot)_{{\mathcal{U}}}$. On ${\mathbb{F}}$, consider a zero-mean
Gaussian family
\[
{\dot{W}}=\left\{  {\dot{W}}(h),\ h\in{\mathcal{U}}\right\}
\]
so that
\[
{\mathbb{E}}\left(  {\dot{W}}(h_{1})\;{\dot{W}}(h_{2})\right)  =(h_{1}%
,h_{2})_{{\mathcal{U}}}.
\]
It suffice, for our purposes, to assume that ${\mathcal{F}}$ is the $\sigma
$-algebra generated by ${\dot{W}}$. Given a real separable Hilbert space $X$,
we denote by $L_{2}({\mathbb{F}};X)$ the Hilbert space of square-integrable
${\mathcal{F}}$-measurable $X$-valued random elements $f$. In particular,
\[
(f,g)_{L_{2}({\mathbb{F}};X)}^{2}:={\mathbb{E}}(f,g)_{X}^{2}.
\]
When $X={\mathbb{R}}$, we write $L_{2}({\mathbb{F}})$ instead of
$L_{2}({\mathbb{F}};{\mathbb{R}})$.

\begin{definition}
\label{def:well} A formal series
\begin{equation}
\dot{W}=\sum_{k}\dot{W}({\mathfrak{u}}_{k})\, {\mathfrak{u}}_{k} ,
\label{eq:wn}%
\end{equation}
where $\left\{  {\mathfrak{u}}_{k},k\geq1\right\}  $ is a complete orthonormal
basis in $\mathcal{U}$,\ is called (Gaussian) \textbf{white noise} on
$\mathcal{U}.$
\end{definition}

The white noise on $\mathcal{U}=L_{2}\left(  G\right)  ,$ where $G$ is a
domain in $\mathbb{R}^{d},$\ is usually referred to as a \textbf{spatial or
space white} \textbf{noise} (on $L_{2}\left(  G\right)  ).$ The space white
noise is of central importance for this paper.

Below, we will introduce a class of spaces that are convenient for treating
nonlinear functionals of white noise, in particular, solutions of SPDEs driven
by white noise.

Given an orthonormal basis ${\mathfrak{U}}=\{{\mathfrak{u}}_{k},k\geq1\}$ in
${\mathcal{U}}$, define a collection $\{\xi_{k},k\geq1\}$ of independent
standard Gaussian random variables so that $\xi_{k}={\dot{W}}({\mathfrak{u}%
}_{k})$. Denote by ${\mathcal{J}}$ the collection of multi-indices $\alpha$
with $\alpha=(\alpha_{1},\alpha_{2},\ldots)$ so that each $\alpha_{k}$ is a
non-negative integer and $|\alpha|:=\sum_{k\geq1}\alpha_{k}<\infty$. For
$\alpha,\beta\in{\mathcal{J}}$, we define
\[
\alpha+\beta=(\alpha_{1}+\beta_{1},\alpha_{2}+\beta_{2},\ldots),\quad
\alpha!=\prod_{k\geq1}\alpha_{k}!.
\]
By $(0)$ we denote the multi-index with all zeroes. By $\varepsilon_{i}$ we
denote the multi-index $\alpha$ with $\alpha_{i}=1$ and $\alpha_{j}=0$ for
$j\not =i$. With this notation, $n\varepsilon_{i}$ is the multi-index $\alpha$
with $\alpha_{i}=n$ and $\alpha_{j}=0$ for $j\not =i$. The following
inequality holds (see Appendix for the proof):
\begin{equation}
\label{HOUZ-ineq}
|\alpha|! \leq \alpha! (2\mathbb{N})^{2\alpha}, \ \ {\rm where}\ \
 (2\mathbb{N})^{2\alpha}=\prod_{k\geq 1} (2k)^{2\alpha_k}.
\end{equation}

Define the collection of random variables $\Xi=\{\xi_{\alpha}, \alpha
\in{\mathcal{J}}\}$ as follows:
\begin{equation}
\label{eq:basis}\xi_{\alpha} = \prod_{k} \left(  \frac{H_{\alpha_{k}}(\xi
_{k})}{\sqrt{\alpha_{k}!}} \right)  ,
\end{equation}
where
\begin{equation}
\label{eq:hermite}H_{n}(x) = (-1)^{n} e^{x^{2}/2}\frac{d^{n}}{dx^{n}}%
e^{-x^{2}/2}%
\end{equation}
is Hermite polynomial of order $n$.

\begin{theorem}[Cameron and Martin \cite{CM}]
\label{th:CM}  The collection $\displaystyle\Xi
=\{\xi_{\alpha},\ \alpha\ \in{\mathcal{J}}\}$ is an orthonormal basis in
$\displaystyle L_{2}({\mathbb{F}})$: if $\displaystyle\eta\in L_{2}%
({\mathbb{F}})$ and $\displaystyle\eta_{\alpha}={\mathbb{E}}(\eta\xi_{\alpha
})$, then $\displaystyle\eta=\sum_{\alpha\in{\mathcal{J}}}\eta_{\alpha}%
\xi_{\alpha}$ and ${\mathbb{E}}|\eta|^{2}=\sum_{\alpha\in{\mathcal{J}}}%
\eta_{\alpha}^{2}.$
\end{theorem}

Expansions with respect to the Cameron-Martin basis $\Xi$ is usually referred
to as Wiener Chaos. Next, we introduce a modification of the Wiener Chaos
expansion which will be called \textit{weighted Wiener Chaos.}

Let ${\mathcal{R}}$ be a bounded linear operator on $L_{2}({\mathbb{F}})$
defined by ${\mathcal{R}}\xi_{\alpha}=r_{\alpha}\xi_{\alpha}$ for every
$\alpha\in{\mathcal{J}}$, where the \emph{weights} $\{r_{\alpha},\ \alpha
\in{\mathcal{J}}\}$ are positive numbers. By Theorem \ref{th:CM},
${\mathcal{R}}$ is bounded if and only if the weights $r_{\alpha}$ are
uniformly bounded from above: $r_{\alpha}<C$ for all $\alpha\in{\mathcal{J}}$,
with $C$ independent of $\alpha$. The inverse operator ${\mathcal{R}}^{-1}$ is
defined by ${\mathcal{R}}^{-1}\xi_{\alpha}=r_{\alpha}^{-1}\xi_{\alpha}$.

We extend ${\mathcal{R}}$ to an operator on $L_{2}({\mathbb{F}};X)$ by
defining ${\mathcal{R}}f$ as the unique element of $L_{2}({\mathbb{F}};X)$ so
that, for all $g\in L_{2}(\mathbb{F};X)$,
\[
{\mathbb{E}}({\mathcal{R}}f,g)_{X}= \sum_{\alpha\in{\mathcal{J}}}r_{\alpha
}{\mathbb{E}}\big((f,g)_{X}\xi_{\alpha}\big).
\]
Denote by ${\mathcal{R}}L_{2}({\mathbb{F}};X)$ the closure of $L_{2}%
({\mathbb{F}};X)$ with respect to the norm
\[
\Vert f\Vert_{{\mathcal{R}}L_{2}({\mathbb{F}};X)}^{2}:=\Vert{\mathcal{R}%
}f\Vert_{L_{2}({\mathbb{F}};X)}^{2}.
\]
Then the elements of ${\mathcal{R}}L_{2}({\mathbb{F}};X)$ can be identified
with a formal series $\sum_{\alpha\in{\mathcal{J}}}f_{\alpha}\xi_{\alpha},$
where $f_{\alpha}\in X$ and $\sum_{\alpha\in{\mathcal{J}}} \Vert f_{\alpha
}\Vert_{X}^{2}r_{\alpha}^{2}<\infty$.

We define the space ${\mathcal{R}}^{-1}L_{2}({\mathbb{F}};X)$ as the dual of
${\mathcal{R}}L_{2}({\mathbb{F}};X)$ relative to the inner product in the
space $L_{2}({\mathbb{R}};X):$
\[
{\mathcal{R}}^{-1}L_{2}({\mathbb{F}};X)=\left\{  g\in L_{2}({\mathbb{F}%
};X):{\mathcal{R}}^{-1}g\in L_{2}({\mathbb{F}};X)\right\}  .
\]
For $f\in{\mathcal{R}}L_{2}({\mathbb{F}};X)$ and $g\in{\mathcal{R}}^{-1}%
L_{2}({\mathbb{F}})$ we define the scalar product%

\begin{equation}
\langle\!\langle{f},{g}\rangle\!\rangle:={\mathbb{E}}\big(({\mathcal{R}%
}f)({\mathcal{R}}^{-1}g)\big)\in X. \label{eq:dualE}%
\end{equation}
In what follows, we will identify the operator ${\mathcal{R}}$ with the
corresponding collection $(r_{\alpha},\ \alpha\in{\mathcal{J}})$. Note that if
$u\in{\mathcal{R}}_{1}L_{2}({\mathbb{F}};X)$ and $v\in{\mathcal{R}}_{2}%
L_{2}({\mathbb{F}};X)$, then both $u$ and $v$ belong to ${\mathcal{R}}%
L_{2}({\mathbb{F}};X)$, where $r_{\alpha}=\min(r_{1,\alpha},r_{2,\alpha})$. As
usual, the argument $X$ will be omitted if $X={\mathbb{R}}$.

Important particular cases of the space ${\mathcal{R}}L_{2}({\mathbb{F}};X)$
correspond to the following weights:

\begin{enumerate}
\item
\[
r_{\alpha}^{2}= \prod_{k=1}^{\infty}q_{k}^{\alpha_{k}},
\]
where $\left\{  q_{k},\, k \geq1\right\}  $ is a non-increasing sequence of
positive numbers with $q_{1}\leq1$ (see \cite{LR_AP,NR});

\item
\begin{equation}
r_{\alpha}^{2}=(\alpha!)^{\rho}(2{\mathbb{N}})^{\ell\alpha},\ \rho\leq
0,\ \ell\leq0,\ \ \mathrm{where\ }(2{\mathbb{N}})^{\ell\alpha}=\prod_{k\geq
1}(2k)^{\ell\alpha_{k}}. \label{eq:S-weight}%
\end{equation}
This set of weights defines Kondratiev's spaces $({\mathcal{S}})_{\rho,\ell
}(X)$.
\end{enumerate}

Now we will sketch the basics of Malliavin calculus on ${\mathcal{R}}%
L_{2}({\mathbb{F}};X).$

Denote by ${\mathbf{D}}$ the \emph{Malliavin derivative} on $L_{2}%
({\mathbb{F}})$ (see e.g. \cite{Nualart}). In particular, if $F:{\mathbb{R}%
}^{N}\rightarrow{\mathbb{R}}$ is a smooth function and $h_{i}\in{\mathcal{U}%
},\ i=1,\ldots N$, then
\begin{equation}
{\mathbf{D}}F({\dot{W}}(h_{1}),\ldots{\dot{W}}(h_{N}))=\sum_{i=1}^{N}%
\frac{\partial F}{\partial x_{i}}({\dot{W}}(h_{1}),\ldots,{\dot{W}}%
(h_{N}))h_{i}\in L_{2}({\mathbb{F}};{\mathcal{U}}). \label{eq:MD}%
\end{equation}
It is known \cite{Nualart} that the domain ${\mathbb{D}}^{1,2}({\mathbb{F}})$
of the operator ${\mathbf{D}}$ is a dense linear subspace of $L_{2}%
({\mathbb{F}})$.

The adjoint of the Malliavin derivative on $L_{2}(\mathbb{F})$ is the
It\^{o}-Skorokhod integral and is traditionally denoted by $\delta$
\cite{Nualart}. We will keep this notation for the extension of this operator
to ${\mathcal{R}}L_{2}({\mathbb{F}};X\otimes{\mathcal{U}})$.

For $f\in{\mathcal{R}}L_{2}({\mathbb{F}};X\otimes{\mathcal{U}})$, we define
${\delta}(f)$ as the unique element of ${\mathcal{R}}L_{2}({\mathbb{F}};X)$
with the property
\begin{equation}
\langle\!\langle{{\delta}(f),\varphi}\rangle\!\rangle={\mathbb{E}%
}({\mathcal{R}}f,{\mathcal{R}}^{-1}{\mathbf{D}}\varphi)_{{\mathcal{U}}}
\label{eq:SkI}%
\end{equation}
for every $\varphi$ satisfying $\varphi\in{\mathcal{R}}^{-1}L_{2}({\mathbb{F}%
})$ and ${\mathbf{D}}\varphi\in{\mathcal{R}}^{-1}L_{2}({\mathbb{F}%
};{\mathcal{U}})$.

Next, we derive the expressions for the Malliavin derivative $\mathbf{D}$ and
its adjoint $\delta$ in the basis $\Xi.$ To begin, we compute ${\mathbf{D}%
}(\xi_{\alpha}).$

\begin{proposition}
\label{prop:Dxi} For each $\alpha\in{\mathcal{J}},$ we have
\begin{equation}
{\mathbf{D}}(\xi_{\alpha})=\sum_{k\geq1}\sqrt{\alpha_{k}}\,\xi_{\alpha
-\varepsilon_{_{k}}}{\mathfrak{u}}_{k}. \label{eq:Dxi}%
\end{equation}

\end{proposition}

\textbf{Proof.} The result follows by direct computation using the property
(\ref{eq:MD}) of the Malliavin derivative and the relation $H_{n}^{\prime
}(x)=nH_{n-1}(x)$ for the Hermite polynomials (cf. \cite{Nualart}).
\endproof

Obviously, the set $\mathcal{J}$ is not invariant with respect to
substraction. In particular, the expression $\alpha-\varepsilon_{k}$ is
undefined if $\alpha_{k}=0$. In (\ref{eq:Dxi}) and everywhere below in this
paper where undefined expressions of this type appear, we use the following
convention: if $\alpha_{k}=0,$ then $\sqrt{\alpha_{k}}\,\xi_{\alpha
-\varepsilon_{k}}=0$.

\begin{proposition}
\label{prop:Ixi} For $\xi_{\alpha}\in\Xi$, $h\in X$, and ${\mathfrak{u}}%
_{k}\in{\mathfrak{U}}$, we have
\begin{equation}
{\delta}(\xi_{\alpha}\, h\otimes{\mathfrak{u}}_{k}) =h\,\sqrt{\alpha_{k}%
+1}\,\xi_{\alpha+\varepsilon_{k}}. \label{eq:Ixi}%
\end{equation}

\end{proposition}

\textbf{Proof.} It is enough to verify (\ref{eq:SkI}) with $f=h\otimes
{\mathfrak{u}}_{k}\,\xi_{\alpha}$ and $\varphi=\xi_{\beta}$, where $h\in X$.
By (\ref{eq:Dxi}),
\[
{\mathbb{E}}(f,{\mathbf{D}}\varphi)_{{\mathcal{U}}}= \sqrt{\beta_{k}}\, h\,
{\mathbb{E}}(\xi_{\alpha}\xi_{\beta-\varepsilon_{k}}) =
\begin{cases}
\sqrt{\alpha_{k} +1}\, h, & \mathrm{if}\ \alpha=\beta-\varepsilon_{k},\\
0, & \mathrm{if } \ \alpha\not =\beta-\varepsilon_{k}.
\end{cases}
\]
In other words,
\[
{\mathbb{E}}(\xi_{\alpha}\, h\otimes{\mathfrak{u}}_{k},{\mathbf{D}}\xi_{\beta
})_{{\mathcal{U}}}=h\,{\mathbb{E}}(\sqrt{\alpha_{k}+1}\xi_{\alpha
+\varepsilon_{k}}\xi_{\beta})
\]
for all $\beta\in{\mathcal{J}}$. \endproof

\begin{remark}
The operator $\delta\mathbf{D}$ is linear and unbounded on $L_{2}(\mathbb{F}%
)$; it follows from Propositions \ref{prop:Dxi} and \ref{prop:Ixi} that the
random variables $\xi_{\alpha}$ are eigenfunctions of this operator:
\begin{equation}
\label{eq:eig-xi}\delta(\mathbf{D}(\xi_{\alpha})) =|\alpha|\xi_{\alpha}.
\end{equation}

\end{remark}

To give an alternative characterization of the operator ${\delta}$, we define
a new operation on the elements of $\Xi$.

\begin{definition}
\label{def:WP} For $\xi_{\alpha}$, $\xi_{\beta}$ from $\Xi$, define the Wick
product
\begin{equation}
\label{eq:WP}\xi_{\alpha}\diamond\xi_{\beta}:=\sqrt{\left(  \frac
{(\alpha+\beta)!}{\alpha!\beta!} \right)  } \xi_{\alpha+\beta}.
\end{equation}

\end{definition}

In particular, taking in (\ref{def:WP}) $\alpha=k\varepsilon_{i}$ and
$\beta=n\varepsilon_{i}$, and using (\ref{eq:basis}), we get
\begin{equation}
\label{eq:WP-hp}H_{k}(\xi_{i})\diamond H_{n}(\xi_{i})=H_{k+n}(\xi_{i}).
\end{equation}

By linearity, we define the Wick product $f\diamond\eta$ for $f\in
{\mathcal{R}} L_{2}({\mathbb{F}};X)$ and $\eta\in{\mathcal{R}} L_{2}%
({\mathbb{F}})$: if $f=\sum_{\alpha\in{\mathcal{J}}}f_{\alpha}\xi_{\alpha}$,
$f_{\alpha}\in X$, and $\eta=\sum_{\alpha\in{\mathcal{J}}}\eta_{\alpha}%
\xi_{\alpha}$, $\eta_{\alpha}\in{\mathbb{R}}$, then
\[
f\diamond\eta= \sum_{\alpha, \beta} f_{\alpha} \eta_{\beta} \xi_{\alpha
}\diamond\xi_{\beta}.
\]

\begin{proposition}
\label{prop:WPgen} If $f\in{\mathcal{R}}L_{2}({\mathbb{F}};X)$ and $\eta
\in{\mathcal{R}}L_{2}({\mathbb{F}})$, then $f\diamond\eta$ is an element of
$\bar{{\mathcal{R}}}L_{2}({\mathbb{F}};X)$ for a suitable operator
$\bar{{\mathcal{R}}}$.
\end{proposition}

\textbf{Proof.} It follows from (\ref{def:WP}) that $f\diamond\eta
=\sum_{\alpha\in{\mathcal{J}}} F_{\alpha}\xi_{\alpha}$ and
\[
F_{\alpha}=\sum_{\beta,\gamma\in{\mathcal{J}}:\beta+\gamma=\alpha}
\sqrt{\left(  \frac{(\beta+\gamma)!}{\beta!\gamma!} \right)  }f_{\beta}%
\eta_{\gamma}.
\]
Therefore, each $F_{\alpha} X$ is an element of $X$, because, for every
$\alpha\in{\mathcal{J}}$, there are only finitely many multi-indices $\beta,
\gamma$ satisfying $\beta+\gamma=\alpha$. It is known \cite[Proposition
7.1]{LR_shir} that
\begin{equation}
\label{eq:HOUZ1}\sum_{\alpha\in{\mathcal{J}}}(2{\mathbb{N}})^{q\alpha}%
<\infty\ \ \mathrm{if\ and\ only\ if\ } q<-1.
\end{equation}
Therefore, $f\diamond\eta\in\bar{{\mathcal{R}}} L_{2}({\mathbb{F}};X)$, where
the operator $\bar{{\mathcal{R}}}$ can be defined using the weights $\bar
{r}_{\alpha}^{2}=(2{\mathbb{N}})^{-2\alpha}/(1+\|F_{\alpha}\|_{X}^{2})$.
\endproof

An immediate consequence of Proposition \ref{prop:Ixi} and Definition
\ref{def:WP} is the following identity:
\begin{equation}
\label{eq:IWP}{\delta}(\xi_{\alpha}h\otimes{\mathfrak{u}}_{k})=h\xi_{\alpha
}\diamond\xi_{k}, \ h\in X.
\end{equation}

Below we summarize the properties of the operator ${\delta}$.

\begin{theorem}
\label{th:dual} If $f$ is an element of ${\mathcal{R}}L_{2}({\mathbb{F}%
};X\otimes{\mathcal{U}})$ so that $f=\sum_{k\geq1}f_{k}\otimes{\mathfrak{u}%
}_{k}$, with $f_{k}=\sum_{\alpha\in{\mathcal{J}}}f_{k,\alpha}\,\xi_{\alpha}%
\in{\mathcal{R}}L_{2}({\mathbb{F}};X)$, then
\begin{equation}
{\delta}(f)=\sum_{k\geq1}f_{k}\diamond\xi_{k}, \label{eq:IWPG}%
\end{equation}
and
\begin{equation}
({\delta}(f))_{\alpha}=\sum_{k\geq1}\sqrt{\alpha_{k}}f_{k,\alpha
-\varepsilon_{k}}. \label{eq:IBasG}%
\end{equation}

\end{theorem}

\textbf{Proof.} By linearity and (\ref{eq:IWP}),
\[
{\delta}(f)=\sum_{k\geq1}\sum_{\alpha\in{\mathcal{J}}} {\delta}(\xi_{\alpha}
f_{k,\alpha}\otimes{\mathfrak{u}}_{k})= \sum_{k\geq1} \sum_{\alpha
\in{\mathcal{J}}} f_{k,\alpha}\xi_{\alpha}\diamond\xi_{k}= \sum_{k\geq1}
f_{k}\diamond\xi_{k},
\]
which is (\ref{eq:IWPG}). On the other hand, by (\ref{eq:Ixi}),
\[
{\delta}(f)=\sum_{k\geq1}\sum_{\alpha\in{\mathcal{J}}} f_{k,\alpha}%
\sqrt{\alpha_{k}+1}\, \xi_{\alpha+\varepsilon_{k}}= \sum_{k\geq1}\sum
_{\alpha\in{\mathcal{J}}} f_{k,\alpha-\varepsilon_{k}}\sqrt{\alpha_{k}}\,
\xi_{\alpha},
\]
and (\ref{eq:IBasG}) follows. \endproof

\begin{remark}
It is not difficult to show that the operator ${\delta}$ can be considered as
an extension of the Skorokhod integral to the weighted spaces ${\mathcal{R}%
}L_{2}({\mathbb{F}};X\otimes{\mathcal{U}})$.
\end{remark}

One way to describe a multi-index $\alpha$ with $|\alpha|=n>0$ is by its
characteristic set $K_{\alpha}$, that is, an ordered $n$-tuple $K_{\alpha
}=\{k_{1},\ldots,k_{n}\}$, where $k_{1}\leq k_{2}\leq\ldots\leq k_{n}$
characterize the locations and the values of the non-zero elements of $\alpha
$. More precisely, $k_{1}$ is the index of the first non-zero element of
$\alpha,$ followed by $\max\left(  0,\alpha_{k_{1}}-1\right)  $ of entries
with the same value. The next entry after that is the index of the second
non-zero element of $\alpha$, followed by $\max\left(  0,\alpha_{k_{2}%
}-1\right)  $ of entries with the same value, and so on.\ For example, if
$n=7$ and $\alpha=(1,0,2,0,0,1,0,3,0,\ldots)$, then the non-zero elements of
$\alpha$ are $\alpha_{1}=1$, $\alpha_{3}=2$, $\alpha_{6}=1$, $\alpha_{8}=3$.
As a result, $K_{\alpha}=\{1,3,3,6,8,8,8\}$, that is, $k_{1}=1,\,k_{2}%
=k_{3}=3,\,k_{4}=6, k_{5}=k_{6}=k_{7}=8$.

Using the notion of the characteristic set, we now state the following analog
of the well-known result of It\^{o} \cite{Ito} connecting multiple Wiener
integrals and Hermite polynomials.

\begin{proposition}
\label{rem:xial} Let $\alpha\in{\mathcal{J}}$ be a multi-index with
$|\alpha|=n\geq1$ and characteristic set $K_{\alpha}=\{k_{1},\ldots, k_{n}\}$.
Then
\begin{equation}
\label{xial-alt}\xi_{\alpha}= \frac{\xi_{k_{1}}\diamond\xi_{k_{2}}%
\diamond\cdots\diamond\xi_{k_{n}}}{\sqrt{\alpha!}}.
\end{equation}

\end{proposition}

\textbf{Proof.} This follows from (\ref{eq:basis}) and (\ref{eq:WP-hp}),
because by (\ref{eq:WP-hp}), for every $i$ and $k$,
\[
H_{k}(\xi_{i})=\underset{k\,\mathrm{times}}{\underbrace{\xi_{i}\diamond
\cdots\diamond\xi_{i}}}.
\]
\endproof

\section{ Evolution Equations Driven by White Noise}

\label{sec3} \setcounter{equation}{0} \setcounter{theorem}{0}

\subsection{The setting}

In this section we study anticipating solutions of stochastic evolution
equations driven by Gaussian white noise on a Hilbert space $\mathcal{U}$.

\begin{definition}
\label{def:Ntrip} The triple $(V,H,V^{\prime})$ of Hilbert spaces is called
\textbf{normal} if and only if

\begin{enumerate}
\item $V\hookrightarrow H \hookrightarrow V^{\prime}$ and both embeddings
$V\hookrightarrow H $ and $H \hookrightarrow V^{\prime}$ are dense and continuous;

\item The space $V^{\prime}$ is the dual of $V$ relative to the inner product
in $H$;

\item There exists a constant $C>0$ so that $|(h, v)_{H}| \leq C\|v\|_{V}%
\|h\|_{V^{\prime}}$ for all $v\in V$ and $h\in H$.
\end{enumerate}
\end{definition}

For example, the Sobolev spaces $(H^{\ell+ \gamma}_{2}({\mathbb{R}}^{d}),
H^{\ell}_{2}({\mathbb{R}}^{d}), H^{\ell- \gamma}_{2}({\mathbb{R}}^{d}))$,
$\gamma>0$, $\ell\in{\mathbb{R}}$, form a normal triple.

Denote by $\langle v^{\prime}, v\rangle$, $v^{\prime}\in V^{\prime}$, $v\in
V$, the duality between $V$ and $V^{\prime}$ relative to the inner product in
$H$. The properties of the normal triple imply that $|\langle v^{\prime},
v\rangle|\leq C\|v\|_{V}\|v^{\prime}\|_{V^{\prime}}$, and, if $v^{\prime}\in
H$ and $v\in V$, then $\langle v^{\prime}, v\rangle= (v^{\prime},v)_{H}.$

We will also use the following notation:
\begin{equation}
{\mathcal{V}}=L_{2}((0,T);V),\ {\mathcal{H}}=L_{2}((0,T);H),\ {\mathcal{V}%
}^{\prime}=L_{2}((0,T);V^{\prime}). \label{sp-notation}%
\end{equation}
Given a normal triple $(V,H,V^{\prime})$, let ${\mathbf{A}}:V\rightarrow
V^{\prime}$ and ${\mathbf{M}}:V\rightarrow V^{\prime}\otimes{\mathcal{U}}$ be
bounded linear operators.

\begin{definition}
\label{def:parab} The solution of the stochastic evolution equation
\begin{equation}
\dot{u}={\mathbf{A}}u+f+{\delta}({\mathbf{M}}u),\ 0<t\leq T, \label{eq:evol}%
\end{equation}
with $f\in{\mathcal{R}}L_{2}({\mathbb{F}};{\mathcal{V}}^{\prime})$ and
$u|_{t=0}=u_{0}\in{\mathcal{R}}L_{2}({\mathbb{F}};H)$, is a process
$u\in{\mathcal{R}}L_{2}({\mathbb{F}};{\mathcal{V}})$ so that, for every
$\varphi$ satisfying $\varphi\in{\mathcal{R}}^{-1}L_{2}({\mathbb{F}})$ and
${\mathbf{D}}\varphi\in{\mathcal{R}}^{-1}L_{2}({\mathbb{F}};{\mathcal{U}})$,
the equality
\begin{equation}
\langle\!\langle{u(t)},{\varphi}\rangle\!\rangle=\langle\!\langle{u_{0}%
},{\varphi}\rangle\!\rangle+\int_{0}^{t}\!\langle\! \langle{{\mathbf{A}%
}u(s)+f(s)+{\delta}({\mathbf{M}}u)(s)},{\varphi}\rangle\!\rangle ds\!
\label{eq:evol-ds}%
\end{equation}
holds in ${\mathcal{V}}^{\prime}$; see (\ref{eq:dualE}) for the definition of
$\langle\!\langle\cdot, \cdot, \rangle\!\rangle$.
\end{definition}

\begin{remark}
\label{re:cont} (a) The solutions described by Definitions \ref{def:parab} and
\ref{def:ell} belong to the class of \textquotedblleft variational solutions",
which is quite typical for partial differential equations (see
\cite{KR,Lions,LM,Roz}, etc.)

(b) Since $\langle\!\langle{u(t)},{\varphi}\rangle\!\rangle\in\mathcal{V}$ and
$\langle\!\langle{u(t)},{\varphi}\rangle\!\rangle_{t}\in\mathcal{V}^{\prime},$
by the standard embedding theorem \textrm{(}see e.g. \cite[Section 1.2.2]%
{LM}\textrm{)} there exists a version of $\langle\!\langle{u(t)},{\varphi
}\rangle\!\rangle\in{\mathbf{C}}\left(  [0,T];H\right)  $. Clearly, one could
also select a version of ${u(t)}$ such that $\langle\!\langle{u(t)},{\varphi
}\rangle\!\rangle\in{\mathbf{C}}\left(  [0,T];H\right)  .$ In the future, we
will consider only this version of the solution. By doing this we ensure that
formula (\ref{eq:evol-ds}) which is understood as an equality in
$\mathcal{V}^{\prime}$ yields $u|_{t=0}=u_{0}\in{\mathcal{R}}L_{2}%
({\mathbb{F}};H).$
\end{remark}

\begin{remark}
To simplify the notations and the overall presentation, we assume that
$\mathbf{A}$ and $\mathbf{M}$ do not depend on time, even though many of the
results in this paper can easily be extended to time-dependent operators.
\end{remark}

Fix an orthonormal basis $\mathfrak{U}$ in ${\mathcal{U}}$. Then, for every
$v\in V$, there exists a collection $v_{k}\in V^{\prime},\ k\geq1$, so that
\[
{\mathbf{M}}v=\sum_{k\geq1}v_{k}\otimes\mathfrak{u}_{k}.
\]
We therefore define the operators ${\mathbf{M}}_{k}:\,V\rightarrow V^{\prime}$
by setting ${\mathbf{M}}_{k}v=v_{k}$ and write
\[
{\mathbf{M}}v=\sum_{k\geq1}({\mathbf{M}}_{k}v)\otimes\mathfrak{u}_{k}.
\]
By (\ref{eq:IWPG}), equation (\ref{eq:evol}) becomes
\begin{equation}
\dot{u}(t)={\mathbf{A}}u(t)+f(t)+{\mathbf{M}}u\left(  t\right)  \diamond
\dot{W}, \label{eq:evol-b}%
\end{equation}
where
\begin{equation}
\mathbf{M}v\diamond\dot{W}:=\sum_{k\geq1}\left(  \mathbf{M}_{k}v\right)
\diamond\xi_{k}. \label{eq:diamondW}%
\end{equation}

\subsection{Equivalence Theorem}

In this section we investigate stochastic Fourier representation of equation
(\ref{eq:evol-b}).

Recall that every process $u=u(t)$ from ${\mathcal{R}}L_{2}({\mathbb{F}%
};{\mathcal{V}})$ is represented by a formal series $u(t)=\sum_{\alpha
\in{\mathcal{J}}}u_{\alpha}(t)\xi_{\alpha}$, with $u_{\alpha}\in{\mathcal{V}}$
and
\begin{equation}
\sum_{\alpha}r_{\alpha}^{2}\left\Vert u_{\alpha}\right\Vert _{\mathcal{V}}%
^{2}<\infty. \label{eq:evolnorm}%
\end{equation}

\begin{theorem}
\label{th: equivalence1}Let $u=\sum_{\alpha\in{\mathcal{J}}}u_{\alpha}%
\xi_{\alpha}$ be an element of ${\mathcal{R}}L_{2}({\mathbb{F}};{\mathcal{V}%
})$. The process $u$ is a solution of equation (\ref{eq:evol}) if and only if
the functions $u_{\alpha}$ have the following properties:

\begin{enumerate}
\item every $u_{\alpha}$ is an element of ${\mathbf{C}} \left(
[0,T];H)\right)  $

\item the system of equalities
\begin{equation}
u_{\alpha}(t)=u_{0,\alpha}+\int_{0}^{t}\left(  \ {\mathbf{A}}u_{\alpha
}(s)+f_{\alpha}(s)+\sum_{k\geq1}\sqrt{\alpha_{k}}\,{\mathbf{M}}_{k}%
u_{\alpha-\varepsilon_{k}}(s)\right)  ds \label{eq:evol-S}%
\end{equation}
holds in $V^{\prime}$ for all $t\in[0,T]$ and $\alpha\in{\mathcal{J}}$.
\end{enumerate}
\end{theorem}

\begin{proof}
Let $u$ be a solution of (\ref{eq:evol}) in ${\mathcal{R}}L_{2}({\mathbb{F}%
};{\mathcal{V}}).$ Taking $\varphi=\xi_{\alpha}$ in (\ref{eq:evol-ds}) and
using relation (\ref{eq:IBasG}), we obtain equation (\ref{eq:evol-S}). By
Remark \ref{re:cont} $u_{\alpha}\in{\mathcal{V}}\bigcap{\mathbf{C}} \left(
[0,T];H)\right)  .$\

Conversely, let $\left\{  u_{\alpha},\alpha\in\mathcal{J}\right\}  $ be a
collection of functions from ${\mathcal{V}}\bigcap{\mathbf{C}}\left(
[0,T];H)\right)  $ satisfying (\ref{eq:evolnorm}) and (\ref{eq:evol-S}). Set
$u\left(  t\right)  :=\sum_{\alpha\in{\mathcal{J}}}u_{\alpha}(t)\xi_{\alpha}$.
$\ $Then, by Theorem \ref{th:dual}, equation (\ref{eq:evol-S}) yields that,
for every $\alpha\in\mathcal{J}$,
\[
\langle\!\langle{u(t)},{\xi_{\alpha}}\rangle\!\rangle=\langle\!\langle{u_{0}%
},{\xi_{\alpha}}\rangle\!\rangle+\int_{0}^{t}\langle\!\langle{{{\mathbf{A}}%
u}\left(  s\right)  {+f}\left(  s\right)  {+{\delta}({\mathbf{M}}u)}(s)}%
,{{\xi}_{\alpha}}\rangle\!\rangle ds.
\]
By continuity, we conclude that for any $\varphi\in{\mathcal{R}}^{-1}%
L_{2}({\mathbb{F}})$ such that ${\mathbf{D}}\varphi\in{\mathcal{R}}^{-1}%
L_{2}({\mathbb{F}};{\mathcal{U}}),$ equality%
\[
\langle\!\langle{u(t)},{\varphi}\rangle\!\rangle=\langle\!\langle{u_{0}%
},{\varphi}\rangle\!\rangle+\int_{0}^{t}\langle\!\langle{{{\mathbf{A}}%
u}\left(  s\right)  {+f}\left(  s\right)  +{\delta}({\mathbf{M}}%
u)(s)},{\varphi}\rangle\!\rangle ds
\]
holds in $\mathcal{V}^{\prime}.$ By Remark \ref{re:cont} $\langle
\!\langle{u(t)},{\varphi}\rangle\!\rangle\in{\mathbf{C}}\left(
[0,T];H\right)  .$
\end{proof}

This simple but very helpful result establishes the equivalence of the
\textquotedblleft physical" (\ref{eq:evol-b}) and the (stochastic) Fourier
(\ref{eq:evol-S}) forms of equation (\ref{eq:evol}). System of equations
(\ref{eq:evol-S}) is often referred in the literature as the \textit{
propagator} of equation (\ref{eq:evol-b}). Note that the propagator is
lower-triangular and can be solved by induction on $|\alpha|$.

\subsection{Existence and uniqueness}

Below, we will present several results on existence and uniqueness of
evolution equations driven by Gaussian white noise.

Before proceeding with general existence-uniqueness problems, we will
introduce two simple examples that indicate the limits of the
\textquotedblleft quality" of solutions of bi-linear SPDEs driven by general
Gaussian white noise.

\begin{example}
\label{ex1} Consider equation
\begin{equation}
u(t)=\phi+\int_{0}^{t}(b\,u(s)\diamond\xi-\lambda u(s))ds, \label{eq:e1}%
\end{equation}
where $\phi,\lambda$ are real numbers, $b$ is a complex number, and $\xi$ is a
standard Gaussian random variable. In other words $\xi$ is Gaussian white
noise on $\mathcal{U}=\mathbb{R}.$\ With only one Gaussian random variable
$\xi$, the set ${\mathcal{J}}$ becomes $\{0,1,2,\ldots\}$ so that
$u(t)=\sum_{n\geq0}u_{(n)}(t)H_{n}(\xi)/\sqrt{n!}$, where $H_{n}$ is Hermite
polynomial of order $n$ (\ref{eq:hermite}). According to (\ref{eq:evol-S}),
\[
u_{(n)}(t)=\phi I_{(n=0)}-\int_{0}^{t}\lambda u_{(n)}(s)ds+\int_{0}^{t}%
b\sqrt{n}u_{(n-1)}(s)ds.
\]
It follows that $u_{(0)}(t)=\phi e^{-\lambda t}$ and then, by induction,
$\displaystyle u_{(n)}(t)=\phi\frac{\left(  b\,t\right)  ^{n}}{\sqrt{n!}%
}e^{-\lambda t}$. As a result,
\[
u(t)=e^{-\lambda t}\Big(\phi+\sum_{n\geq1}\frac{\left(  b\,t\right)  ^{n}}%
{n!}H_{n}(\xi)\Big)=\phi e^{-\lambda t+(b\,t\xi-\left\vert b\right\vert
^{2}t^{2}/2)}.
\]

\end{example}

Obviously, the solution of the equation is  square integrable on any fixed
time interval. However, as the next example indicates, the solutions of SPDEs
driven by stationary noise are much more intricate then the non-anticipating,
or adapted, solutions of SPDEs driven by space-time white noise.

\begin{example}
\label{ex3} With $\xi$ as in the previous examples, consider a partial
differential equation
\begin{equation}
u_{t}(t,x)=au_{xx}\left(  t,x\right)  +\left(  \beta u\left(  t,x\right)
+\sigma u_{x}\left(  t,x\right)  \right)  \diamond\xi,\ t>0,\ x\in{\mathbb{R}%
}, \label{eq:ex3}%
\end{equation}
with some initial condition $u_{0}\in L_{2}({\mathbb{R}})$. By taking the
Fourier transform and using the results of Example \ref{ex1} with $\phi
=\hat{u}_{0}(y)$, $\lambda=-ay^{2}$, $b=\beta+\sqrt{-1} y \sigma$, we find
\[%
\begin{split}
\hat{u}_{t}(t,y)  &  =-y^{2}a\hat{u}\left(  t\right)  +\left(  \beta+\sqrt
{-1}y\sigma\right)  \hat{u}\left(  t,y\right)  \diamond\xi;\\
\hat{u}(t,y)  &  =\hat{u}_{0}(y)\exp\left(  -tay^{2}+\left(  \sigma^{2}%
y^{2}-\beta^{2}\right)  t^{2}/2+\sqrt{-1}\beta\sigma yt^{2}+\left(  \sqrt
{-1}\sigma y+\beta\right)  t\xi\right)  .
\end{split}
\]
If $\sigma=0,$ i.e. the \textquotedblleft{}diffusion" operator in equation
(\ref{eq:ex3}) is of order zero, then the solution belongs to $L_{2}%
({\mathbb{F}};L_{2}({\mathbb{R}}))$ for all $t.$ However, if $\sigma>0,$ then
the solution $u(t,\cdot)$ will, in general, belong to $L_{2}({\mathbb{F}%
};L_{2}({\mathbb{R}}))$ only for $t\leq2a/\sigma^{2}$. This blow-up in finite
time is in sharp contrast with the solution of the equation
\begin{equation}
u_{t}=au_{xx}+\sigma u_{x}\diamond\dot{w}, \label{eq:superpar}%
\end{equation}
driven by the standard one-dimensional white noise $\dot{w}\left(  t\right)
=\partial_{t}W\left(  t\right)  $, where $W\left(  t\right)  $ is the
one-dimensional Brownian motion; a more familiar way of writing
(\ref{eq:superpar}) is in the It\^{o} form
\begin{equation}
\label{eq:superparI}du=au_{xx}dt+\sigma u_{x} dW(t).
\end{equation}
It is well known \textrm{(}see, for example, \cite{Roz}\textrm{)} that the
solution of (\ref{eq:superparI}) belongs to $L_{2}({\mathbb{F}};L_{2}%
({\mathbb{R}}))$ for every $t>0$ as long as $u_{0}\in L_{2}({\mathbb{R}})$
and
\begin{equation}
a-\sigma^{2}/2\geq0. \label{co:superpar}%
\end{equation}

\end{example}

The existence of a square integrable (global) solution of an It\^{o}'s SPDE
with square integrable initial condition hinges on the parabolic condition
which in the case of equation (\ref{eq:superpar}) is given by
(\ref{co:superpar}). Example \ref{ex3} shows that this condition is not in any
way sufficient for SPDEs involving a Skorokhod-type integral. The next theorem
provides sufficient conditions for the existence and uniqueness of a solution
to equation (\ref{eq:evol-b}) in the space ${\mathcal{R}}L_{2}({\mathbb{F}%
};{\mathcal{V}}),$ which appears to be a reasonable extension of the class of
square integrable solutions.

Firstly, we introduce an additional assumption on the operator ${\mathbf{A}}$
that will be used throughout this section:

(A):\textbf{ }\textit{For every }$U_{0}\in H$\textit{ and }$F\in
\mathcal{V}^{\prime}:=L_{2}((0,T);V^{\prime})$\textit{, there exists a
function }$U\in\mathcal{V}$\textit{ that solve the deterministic equation }%
\begin{equation}
\partial_{t}U(t)={\mathbf{A}}U(t)+F(t),\ U\left(  0\right)  =U_{0},
\label{eq: determ}%
\end{equation}
\textit{ and there exists a constant }$C=C\left(  {\mathbf{A}},T\right)
$\textit{ so that }%
\begin{equation}
\Vert U\Vert_{{\mathcal{V}}}\leq C({\mathbf{A}},T)\big(\Vert U_{0}\Vert
_{H}+\Vert F\Vert_{{\mathcal{V}}^{\prime}}\big). \label{det-reg}%
\end{equation}

\begin{remark}
\label{re: cont1} Assumption (A) implies that a solution of equation
(\ref{eq: determ}) is unique and belongs to ${\mathbf{C}} \left(
(0,T);H\right)  $ (cf. Remark \ref{re:cont}). The assumption also implies that
the operator ${\mathbf{A}}$ generates a semi-group $\Phi=\Phi_{t}, \, t\geq0,$
and, for every $v\in{\mathcal{V}}$,
\begin{equation}
\int_{0}^{T}\left\Vert \int_{0}^{t}\ \Phi_{t-s}{\mathbf{M}}_{k}v\left(
s\right)  ds\right\Vert _{V}^{2}dt\leq C_{k}^{2}\left\Vert v\right\Vert
_{\mathcal{V}}^{2}, \label{Ck-def}%
\end{equation}
with numbers $C_{k}$ independent of $v$.
\end{remark}

\begin{remark}
\label{rm:evol} There are various types of assumptions on the operator
${\mathbf{A}}$ that yield the statement of the assumption (A). In particular,
(A) holds if the operator ${\mathbf{A}}$ \ is coercive in $\left(
V,H,V^{\prime}\right)  $:
\[
\langle{\mathbf{A}}v,v\rangle+\gamma\Vert v\Vert_{V}^{2}\leq C\Vert v\Vert
_{H}^{2}%
\]
for every $v\in V$, where $\gamma>0$ and $C\in{\mathbb{R}}$ are both
independent of $v$.
\end{remark}

\begin{theorem}
\label{th:evol} Assume(A). Consider equation (\ref{eq:evol-b}) in which
$u_{0}\in\bar{{\mathcal{R}}}L_{2}({\mathbb{F}};H)$, $f\in\bar{{\mathcal{R}}%
}L_{2}({\mathbb{F}};{\mathcal{V}}^{\prime})$ for some operator $\bar
{{\mathcal{R}}},$ and each ${\mathbf{M}}_{k}$ is a bounded linear operator
from $V$ to $V^{\prime}$.

Then there exist an operator ${\mathcal{R}}$ and a unique solution
$u\in{\mathcal{R}}L_{2}({\mathbb{F}};{\mathcal{V}})$ of (\ref{eq:evol-b}).
\end{theorem}

\textbf{Proof.} By Theorem \ref{th: equivalence1}, it suffices to prove that
the propagator (\ref{eq:evol-S}) has a unique solution $\left(  u_{\alpha
}\left(  t\right)  \right)  _{\alpha\in\mathcal{J}}$ such that for each
$\alpha,u_{\alpha}\in\mathcal{V}\bigcap{\mathbf{C}} \left(  \left[
0,T\right]  ;H\right)  $ and $u:=\sum_{\alpha\in\mathcal{J}}u_{\alpha}%
\xi_{\alpha}\in{\mathcal{R}}L_{2}({\mathbb{F}};{\mathcal{V}}).$

For $\alpha=(0)$, that is, when $|\alpha|=0$, equation (\ref{eq:evol-S})
reduces to%
\[
u_{\left(  0\right)  }(t)=u_{0,(0)}+\int_{0}^{t}\left(  {\mathbf{A}}u_{\left(
0\right)  }+f_{\left(  0\right)  }\right)  (s)ds.
\]
By (A), this equation has a unique solution and
\[
\Vert u_{(0)}\Vert_{{\mathcal{V}}}\leq C({\mathbf{A}},T)\left(  \Vert
u_{0,(0)}\Vert_{H}+\Vert f_{(0)}\Vert_{{\mathcal{V}}^{\prime}}\right)  .
\]
Using assumption (A), it follows by induction on $|\alpha|$ that, for every
$\alpha\in{\mathcal{J}}$, equation
\begin{equation}
\partial_{t}u_{\alpha}\left(  t\right)  ={\mathbf{A}} u_{\alpha}\left(
t\right)  +f_{\alpha}\left(  t\right)  +\sum_{k\geq1}\sqrt{\alpha_{k}%
}\,{\mathbf{M}}_{k}u_{\alpha-\varepsilon_{k}}\left(  t\right)  ,\text{
}u_{\alpha}\left(  0\right)  =u_{0,\alpha} \label{u-alpha-sol}%
\end{equation}
has a unique solution in $\mathcal{V}\bigcap{\mathbf{C}} \left(  \left[
0,T\right]  ;H\right)  .$ \ Moreover, by (\ref{det-reg}),
\[
\Vert u_{\alpha}\Vert_{{\mathcal{V}}}\leq{\overline{C}}({\mathbf{A}%
},{\mathbf{M}},T)\left(  \Vert u_{0,\alpha}\Vert_{H}+\Vert f_{\alpha}%
\Vert_{{\mathcal{V}}^{\prime}}+\sum_{k\geq1}\sqrt{\alpha_{k}}\Vert
u_{\alpha-\varepsilon_{k}}\Vert_{{\mathcal{V}}}\right)  .
\]
Since only finitely many of $\alpha_{k}$ are different from $0$, we conclude
that $\Vert u_{\alpha}\Vert_{{\mathcal{V}}}<\infty$ for all $\alpha
\in{\mathcal{J}}$.

Define the operator ${\mathcal{R}}$ on $L_{2}({\mathbb{F}})$ using the
weights
\[
r_{\alpha}=\min\left(  \bar{r}_{\alpha}, \frac{(2{\mathbb{N}})^{-\kappa\alpha
}}{1+\Vert u_{\alpha}\Vert_{{\mathcal{V}}}}\right)  ,
\]
where $\kappa>1/2$ (cf. (\ref{eq:S-weight})). Then $u(t):=\sum_{\alpha
\in{\mathcal{J}}}u_{\alpha}(t)\xi_{\alpha}$ is a solution of (\ref{eq:evol-b})
and, by (\ref{eq:HOUZ1}), belongs to ${\mathcal{R}}L_{2}({\mathbb{F}%
};{\mathcal{V}})$. \endproof

While Theorem \ref{th:evol} establishes that under very broad assumptions one
can find an operator ${\mathcal{R}}$ such that equation (\ref{eq:evol-b}) has
a unique solution in ${\mathcal{R}}L_{2}({\mathbb{F}};{\mathcal{V}})$, the
choice of the operator ${\mathcal{R}}$ is not sufficiently explicit (because
of the presence of $\|u_{\alpha}\|_{\mathcal{V}}$) and is not necessarily optimal.

Consider equation (\ref{eq:evol-b}) with non-random $f$ and $u_{0}$. In this
situation, it is possible to find more constructive expression for $r_{\alpha
}$ and to derive explicit formulas, both for ${\mathcal{R}}u$ and for each
individual $u_{\alpha}$.

\begin{theorem}
\label{th:evol-nr} If $u_{0}$ and $f$ are non-random, then the following holds:

\begin{enumerate}
\item the coefficient $u_{\alpha}$, corresponding to the multi-index $\alpha$
with $|\alpha|=n\geq1$ and characteristic set $K_{\alpha}=\{k_{1},\ldots
,k_{n}\},$ is given by
\begin{equation}
u_{\alpha}\left(  t\right)  =\frac{1}{\sqrt{\alpha!}}\sum_{\sigma
\in{\mathcal{P}}_{n}}\int_{0}^{t}\int_{0}^{s_{n}}\ldots\int_{0}^{s_{2}}%
\Phi_{t-s_{n}}{\mathbf{M}}_{k_{\sigma(n)}}\cdots\Phi_{s_{2}-s_{1}}{\mathbf{M}%
}_{k_{\sigma(1)}}u_{(0)}\left(  s_{1}\right)  ds_{1}\ldots ds_{n},
\label{coef-evol}%
\end{equation}
where

\begin{itemize}
\item ${\mathcal{P}}_{n}$ is the permutation group of the set $(1,\ldots, n)$;

\item $\Phi_{t}$ is the semigroup generated by ${\mathbf{A}}$;

\item $u_{(0)}(t)=\Phi_{t}u_{0}+\int_{0}^{t}\Phi_{t-s}f(s)ds$.
\end{itemize}

\item the weights $r_{\alpha}$ can be taken in the form
\begin{equation}
\label{r-alpha-evol}r_{\alpha}=\frac{q^{\alpha}}{\sqrt{|\alpha|!}%
},\ \mathrm{where\ } q^{\alpha}=\prod_{k=1}^{\infty}q_{k}^{\alpha_{k}},
\end{equation}
and the numbers $q_{k}, \, k\geq1,$ are chosen so that $\sum_{k\geq1}
q_{k}^{2}C^{2}_{k}<1$, with $C_{k}$ from (\ref{Ck-def}).

\item With $q_{k}$ and $r_{\alpha}$ from (\ref{r-alpha-evol}), we have
\begin{equation}
\label{multiint-evol1}\sum_{|\alpha|=n} q^{\alpha} u_{\alpha}(t)\xi_{\alpha}=
\int_{0}^{t}\int_{0}^{s_{n}}\!\!\!\!\! \ldots\int_{0}^{s_{2}} \Phi_{t-s_{n}%
}{\delta}({\overline{{\mathbf{M}}}}\Phi_{s_{n}-s_{n-1}} {\delta}(\ldots
{\delta} ({\overline{{\mathbf{M}}}} u_{(0)}))\ldots)ds_{1}\ldots ds_{n-1}
ds_{n},
\end{equation}
where ${\overline{{\mathbf{M}}}} =(q_{1}{\mathbf{M}}_{1}, q_{2}{\mathbf{M}%
}_{2}, \ldots)$, and
\begin{equation}
\label{multiint-evol2}%
\begin{split}
&  {\mathcal{R}} u(t)=u_{(0)}(t)\\
&  +\sum_{n=1}^{\infty}\frac{1}{2^{n}\sqrt{n!}} \int_{0}^{t}\int_{0}^{s_{n}%
}\!\!\!\! \ldots\int_{0}^{s_{2}} \Phi_{t-s_{n}}{\delta}({\overline
{{\mathbf{M}}}}\Phi_{s_{n}-s_{n-1}} {\delta}(\ldots{\delta} ({\overline
{{\mathbf{M}}}} u_{0}(s_{1})))\ldots)ds_{1}\ldots ds_{n-1} ds_{n}.
\end{split}
\end{equation}

\end{enumerate}
\end{theorem}

\textbf{Proof.} If $u_{0}$ and $f$ are deterministic, then equation
(\ref{eq:evol-S}) becomes
\begin{align}
u_{(0)}(t)  &  =u_{0}+\int_{0}^{t}{\mathbf{A}}u_{(0)}(s)ds+\int_{0}%
^{t}f(s)ds,\ |\alpha|=0;\label{eq:evol-S-det0}\\
u_{\alpha}(t)  &  =\int_{0}^{t}{\mathbf{A}}u_{\alpha}(s)ds+\sum_{k\geq1}%
\sqrt{\alpha_{k}}\int_{0}^{t}{\mathbf{M}}_{k}u_{\alpha-\varepsilon_{k}%
}(s)ds,\ |\alpha|>0. \label{eq:evol-S-det1}%
\end{align}
Define ${\widetilde{u}}_{\alpha}=\sqrt{\alpha!}\,u_{\alpha}$. Then
${\widetilde{u}}_{(0)}=u_{(0)}$ and, for $|\alpha|>0$, (\ref{eq:evol-S-det1})
implies
\[
{\widetilde{u}}_{\alpha}(t)=\int_{0}^{t}{\mathbf{A}}{\widetilde{u}}_{\alpha
}(s)ds+\sum_{k\geq1}\int_{0}^{t}{\alpha_{k}}\,{\mathbf{M}}_{k}{\widetilde{u}%
}_{\alpha-\varepsilon_{k}}(s)ds
\]
or
\[
{\widetilde{u}}_{\alpha}(t)=\sum_{k\geq1}\alpha_{k}\int_{0}^{t}\Phi
_{t-s}{\mathbf{M}}_{k}\tilde{u}_{\alpha-\varepsilon_{k}}(s)ds=\sum_{k\in
K_{\alpha}}\int_{0}^{t}\Phi_{t-s}{\mathbf{M}}_{k}\tilde{u}_{\alpha
-\varepsilon_{k}}(s)ds.
\]
By induction on $n$,
\[
{\widetilde{u}}_{\alpha}(t)=\sum_{\sigma\in{\mathcal{P}}_{n}}\int_{0}^{t}%
\int_{0}^{s_{n}}\ldots\int_{0}^{s_{2}}\Phi_{t-s_{n}}{\mathbf{M}}%
_{k_{\sigma(n)}}\cdots\Phi_{s_{2}-s_{1}}{\mathbf{M}}_{k_{\sigma(1)}}%
u_{(0)}ds_{1}\ldots ds_{n},
\]
and (\ref{coef-evol}) follows.

Since (\ref{multiint-evol2}) follows directly from (\ref{multiint-evol1}), it
remains to establish (\ref{multiint-evol1}). To this end, define
\[
U_{n}(t)=\sum_{|\alpha|=n}q^{\alpha}u_{\alpha}(t)\xi_{\alpha},\ n\geq0.
\]
Let us first show that, for each $n\geq1$, $U_{n}\in L_{2}({\mathbb{F}%
};{\mathcal{V}})$. Indeed, for $\alpha=(0)$, $u_{\alpha}(0)=u_{0}$,
$f_{\alpha}=f$ and
\[
u_{(0)}(t)=\Phi_{t}u_{0}+\int_{0}^{t}\Phi_{t-s}f(s)ds.
\]
By (\ref{det-reg}), we have
\begin{equation}
\label{aux3.16}\Vert u_{(0)}\Vert_{{\mathcal{V}}}\leq C({\mathbf{A}},T)\left(
\Vert u_{0}\Vert_{H}+\Vert f\Vert_{{\mathcal{V}}^{\prime}}\right)  .
\end{equation}
When $|\alpha|\geq1$, $f_{\alpha}=0$ and the solution of (\ref{eq:evol-S-det1}%
) is given by
\begin{equation}
u_{\alpha}(t)=\sum_{k\geq1}\sqrt{\alpha_{k}}\int_{0}^{t}\Phi_{t-s}{\mathbf{M}%
}_{k}u_{\alpha-\varepsilon_{k}}(s)ds. \label{eq: u-alpha-sol}%
\end{equation}

By (\ref{coef-evol}), together with (\ref{det-reg}), (\ref{aux3.16}), and
(\ref{Ck-def}), we have
\begin{equation}
\Vert u_{\alpha}\Vert_{{\mathcal{V}}}^{2}\leq C^{2}({\mathbf{A}}%
,T)\frac{(|\alpha|!)^{2}}{\alpha!}\big(\Vert u_{0}\Vert_{H}^{2}+\Vert
f\Vert_{{\mathcal{V}}^{\prime}}^{2}\big)\,\prod_{k\geq1}C_{k}^{2\alpha_{k}%
}.\label{aux-evol-nr000}%
\end{equation}
By the multinomial formula,
\begin{equation}
\label{multinomial}
\left(\sum_{k\geq 1} x_k\right)^n=\sum_{|\alpha|=n}\left(\frac{n!}{\alpha!}
\prod_{k\geq 1} x_k^{\alpha_k}\right).
\end{equation}
Then
\[%
\begin{split}
\sum_{|\alpha|=n}q^{2\alpha}\Vert u_{\alpha}\Vert_{{\mathcal{V}}}^{2} &  \leq
C^{2}({\mathbf{A}},T)\big(\Vert u_{0}\Vert_{H}^{2}+\Vert f\Vert_{{\mathcal{V}}^{\prime}}^{2}\big)\,n!\,\sum_{|\alpha|=n}
\left(\frac{n!}{\alpha!}\prod_{k\geq1}(C_{k}q_{k})^{2\alpha_{k}}\right)\\
&  =C^{2}({\mathbf{A}},T)\big(\Vert u_{0}\Vert_{H}^{2}+\Vert f\Vert_{{\mathcal{V}}^{\prime}}^{2}\big)\,n!\,\left(  \sum_{k\geq1}C_{k}^{2}q_{k}^{2}\right)^{n}<\infty,
\end{split}
\]
because of the selection of $q_{k}$, and so $U_{n}\in L_{2}({\mathbb{F}%
};{\mathcal{V}})$. Moreover, if the weights $r_{\alpha}$ are defined by
(\ref{r-alpha-evol}), then
\[
\sum_{\alpha\in{\mathcal{J}}}r_{\alpha}^{2}\Vert u_{\alpha}\Vert
_{{\mathcal{V}}}^{2}=\sum_{n\geq0}\sum_{|\alpha|=n}r_{\alpha}^{2}\Vert
u_{\alpha}\Vert_{{\mathcal{V}}}^{2}\leq C^{2}({\mathbf{A}},T)\big(\Vert u_{0}\Vert_{H}^{2}+\Vert f\Vert_{{\mathcal{V}}^{\prime}}^{2}\big)\sum_{n\geq1}\left(
\sum_{k\geq1}C_{k}^{2}q_{k}^{2}\right)  ^{n}<\infty
\]
because of the assumption $\sum_{k\geq1}C_{k}^{2}q_{k}^{2}<1$.

Next, the definition of $U_{n}(t)$ and (\ref{eq: u-alpha-sol}) imply that
(\ref{multiint-evol1}) is equivalent to
\begin{equation}
U_{n}(t)=\int_{0}^{t}\Phi_{t-s}{\delta}({\overline{{\mathbf{M}}}}%
U_{n-1}(s))ds,\ n\geq1. \label{evol-aux11}%
\end{equation}
Accordingly, we will prove (\ref{evol-aux11}). For $n=1$, we have
\[
U_{1}(s)=\sum_{k\geq1}q_{k}u_{\varepsilon_{k}}(t)\xi_{k}=\sum_{k\geq1}\int
_{0}^{t}q_{k}\Phi_{t-s}{\mathbf{M}}_{k}u_{(0)}\xi_{k}dt=\int_{0}^{t}\Phi
_{t-s}{\delta}({\overline{{\mathbf{M}}}}U_{0}(s))ds,
\]
where the last equality follows from (\ref{eq:IWPG}). More generally, for
$n>1$ we have by definition of $U_{n}$ that
\[
(U_{n})_{\alpha}(t)=%
\begin{cases}
q^{\alpha}u_{\alpha}(t), & \mathrm{if}\ |\alpha|=n,\\
0, & \mathrm{otherwise}.
\end{cases}
\]
From the equation
\[
q^{\alpha}u_{\alpha}(t)=\int_{0}^{t}{\mathbf{A}}q^{\alpha}u_{\alpha}%
(s)ds+\sum_{k\geq1}\int_{0}^{t}q_{k}\sqrt{\alpha_{k}}\,{\mathbf{M}}%
_{k}q^{\alpha-\varepsilon_{k}}u_{\alpha-\varepsilon_{k}}(s)ds
\]
we find
\[%
\begin{split}
(U_{n}(t))_{\alpha}  &  =%
\begin{cases}
\displaystyle \sum_{k\geq1}\sqrt{\alpha_{k}}q_{k}\int_{0}^{t}\Phi
_{t-s}{\mathbf{M}}_{k}q^{\alpha-\varepsilon_{k}}u_{\alpha-\varepsilon_{k}%
}(s)ds, & \mathrm{if}\ |\alpha|=n,\\
0, & \mathrm{otherwise}.
\end{cases}
\\
&  =\sum_{k\geq1}\sqrt{\alpha_{k}}\,\int_{0}^{t}\Phi_{t-s}{\overline
{{\mathbf{M}}}}_{k}(U_{n-1}(s))_{\alpha-\varepsilon_{k}}ds,
\end{split}
\]
and then (\ref{evol-aux11}) follows from (\ref{eq:IBasG}). Theorem
\ref{th:stat-nr} is proved. \endproof

Formula (\ref{multiint-evol1}) is similar to the multiple Wiener integral
representation of the solution of a stochastic parabolic equation driven by
the Wiener process; see \cite[Theorem 3.8]{LR_AP}.

\begin{example}
Consider the equation
\begin{equation}
u(t,x)=u_{0}(x)+\int_{0}^{t}u_{xx}(s,x)ds+\sum_{k\geq1}\int_{0}^{t}\sigma
_{k}u_{xx}(s,x)\diamond\xi_{k}ds. \label{eq:two-der}%
\end{equation}
With no loss of generality assume that $\sigma_{k}\not =0$ for all $k$.
Standard properties of the heat kernel imply assumption (A) and inequality
(\ref{Ck-def}) with $C_{k}=\sigma_{k}^{2}$. Then the conclusions of Theorem
\ref{th:evol-nr} hold, and we can take $q_{k}^{2}=k^{-2}4^{-k}(1+\sigma
_{k}^{2})^{-k}$. Note that Theorem \ref{th:evol-nr} covers equation
(\ref{eq:two-der}) with no restrictions on the numbers $\sigma_{k}$.
\end{example}

\bigskip In the existing literature on the subject, equations of the type
(\ref{eq:evol-b}) are considered only under the assumption

\textit{(H): each }${\mathbf{M}}_{k}$\textit{ is a bounded linear operators
from }$V$\textit{ to }$H.$\newline Obviously this assumption rules out
equation (\ref{eq:two-der}) but still covers equation (\ref{eq:ex3}). \

Of course, Theorem \ref{th:evol-nr} does not rule out a possibility of a
better-behaving solution under additional assumptions on the operators
${\mathbf{M}}_{k}$. Indeed, it was shown in \cite{LR_shir} that if (H) is
assumed and the space-only Gaussian noise in equation (\ref{eq:evol-b})\ is
replaced by the space-time white noise, then a more delicate analysis of
equation (\ref{eq:evol-b}) is possible. In particular, the solution can belong
to a much smaller Wiener chaos space even if $u_{0}$ and $f$ are not deterministic.

If the operators $\mathbf{M}_{k}$ are bounded in $H$ (see e.g. equation
(\ref{eq:ex3}) with $\sigma=0),$ then, as the following theorem shows, the
solutions can be square integrable (cf. \cite{Hu}).

\begin{theorem}
\label{th:evol-sp} Assume that the operator ${\mathbf{A}}$ satisfies
\begin{equation}
\label{eq:cond(A)}\langle{\mathbf{A}} v,v\rangle+ \kappa\|v\|^{2}_{V}\leq
C_{A}\|v\|_{H}^{2}%
\end{equation}
for every $v\in V$, with $\kappa>0,\; C_{A}\in{\mathbb{R}}$ independent of
$v$, and assume that each ${\mathbf{M}}_{k}$ is a bounded operator on $H$ so
that $\|{\mathbf{M}}_{k}\|_{H\to H}\leq c_{k}$ and
\begin{equation}
\label{eq:cond(M)}C_{M}:=\sum_{k\geq1}c_{k}^{2}<\infty.
\end{equation}
If $f\in{\mathcal{V}}^{\prime}$ and $u_{0}\in H$ are non-random, then there
exists a unique solution $u$ of (\ref{eq:evol-b}) so that $u(t)\in
L_{2}({\mathbb{F}};H)$ for every $t$ and
\begin{equation}
\label{eq:evol-sp1}{\mathbb{E}} \|u(t)\|_{H}^{2}\leq C(C_{A},C_{M}%
,\kappa,t)\left(  \int_{0}^{t}\|f(s)\|_{V^{\prime}}^{2}ds+\|u_{0}\|_{H}^{2}
\right)  .
\end{equation}

\end{theorem}

\textbf{Proof.} Existence and uniqueness of the solution follow from Theorem
\ref{th:evol} and Remark \ref{rm:evol}, and it remains to establish
(\ref{eq:evol-sp1}).

It follows from (\ref{eq:evol-S}) that
\begin{equation}
u_{\alpha}=\frac{1}{\sqrt{\alpha!}}\sum_{k\in K_{\alpha}}^{|\alpha|}\int
_{0}^{t}\Phi_{t-s}{\mathbf{M}}_{k}u_{\alpha-\varepsilon_{k}}(s)ds,
\label{eq:evol-sp2}%
\end{equation}
where $\Phi$ is the semi-group generated by ${\mathbf{A}}$ and $K_{\alpha}$ is
the characteristic set of $\alpha$. Assumption (\ref{eq:cond(A)}) implies that
$\Vert\Phi_{t}\Vert_{H\rightarrow H}\leq e^{pt}$ for some $p\in{\mathbb{R}}$.
A straightforward calculation using relation (\ref{eq:evol-sp2}) and induction
on $|\alpha|$ shows that
\begin{equation}
\Vert u_{\alpha}(t)\Vert_{H}\leq e^{pt}\frac{t^{|\alpha|}c^{\alpha}}%
{\sqrt{\alpha!}}\Vert u_{(0)}\Vert_{H}, \label{eq:evol-sp3}%
\end{equation}
where $c^{\alpha}=\prod_{k}c_{k}^{\alpha_{k}}$ and $u_{(0)}(t)=\Phi_{t}%
u_{0}+\int_{0}^{t}\Phi_{t-s}f(s)ds$. Assumption (\ref{eq:cond(A)}) implies
that $\Vert u_{(0)}\Vert_{H}^{2}\leq C(C_{A},\kappa,t)\left(  \int_{0}%
^{t}\Vert f(s)\Vert_{V^{\prime}}^{2}ds+\Vert u_{0}\Vert_{H}^{2}\right)  $. To
establish (\ref{eq:evol-sp1}), it remans to observe that
\[
\sum_{\alpha\in{\mathcal{J}}}\frac{c^{2\alpha}t^{2|\alpha|}}{\alpha!}%
=e^{C_{M}t^{2}}.
\]
Theorem \ref{th:evol-sp} is proved. \endproof

\begin{remark}
\label{rm:evol-sp} Taking ${\mathbf{M}}_{k}u=c_{k}u$ shows that, in general,
bound (\ref{eq:evol-sp3}) cannot be improved. When condition (\ref{eq:cond(M)}%
) does not hold, a bound similar to (\ref{eq:evol-sp1}) can be established in
a weighted space ${\mathcal{R}}L_{2}({\mathbb{F}};H)$, for example with
$r_{\alpha}=q^{\alpha}$, where $q_{k}=1/(2^{k}(1+c_{k}))$. For special
operators ${\mathbf{M}}_{k}$, a more delicate analysis might be possible; see,
for example, \cite{Hu}.
\end{remark}

If $f$ and $u_{0}$ are not deterministic, then the solution of
(\ref{eq:evol-b}) might not satisfy
\[
{\mathbb{E}} \|u(t)\|_{H}^{2}\leq C(C_{A},C_{M},\kappa,t) \left(  \int_{0}%
^{t}{\mathbb{E}}\|f(s)\|_{V^{\prime}}^{2}ds +{\mathbb{E}}\|u_{0}\|_{H}^{2}
\right)
\]
even if all other conditions of Theorem \ref{th:evol-sp} are fulfilled. An
example can be constructed similar to Example 9.7 in \cite{LR_shir}: an
interested reader can verify that the solution of the equation $u(t)=u_{0}%
+\int_{0}^{t}u(s)\diamond\xi\, ds$, where $\xi$ is a standard Gaussian random
variable and $\displaystyle u_{0}=\sum_{n\geq0} a_{n}\frac{ H_{n}(\xi)}%
{\sqrt{n!}}$, satisfies $\displaystyle {\mathbb{E}} u^{2}(1)\ge\frac{1}{10}
\sum_{n\geq1} a_{n}^{2}e^{\sqrt{n}}$. For equations with random input, one
possibility is to use the spaces $({\mathcal{S}})_{-1,q};$ see
(\ref{eq:S-weight}). Examples of the corresponding results are Theorems
\ref{th:ell-sp2} and \ref{th:wn-evol} below and Theorem 9.8 in \cite{LR_shir}.

\section{Stationary equations}

\label{sec4} \setcounter{equation}{0} \setcounter{theorem}{0}

\subsection{Definitions and Analysis}

The objective of this section is to study stationary stochastic equation
\begin{equation}
{\mathbf{A}}u+{\delta}({\mathbf{M}}u)=f. \label{eq:ell}%
\end{equation}

\begin{definition}
\label{def:ell} The solution of equation (\ref{eq:ell}) with $f\in
{\mathcal{R}}L_{2}({\mathbb{F}};V^{\prime})$, is a random element
$u\in{\mathcal{R}}L_{2}({\mathbb{F}};V)$ so that, for every $\varphi$
satisfying $\varphi\in{\mathcal{R}}^{-1}L_{2}({\mathbb{F}})$ and ${\mathbf{D}%
}\varphi\in{\mathcal{R}}^{-1}L_{2}({\mathbb{F}};{\mathcal{U}})$, the equality
\begin{equation}
\langle\!\langle{{\mathbf{A}}u},{\varphi}\rangle\!\rangle+\langle
\!\langle{{\delta}({\mathbf{M}}u)},{\varphi}\rangle\!\rangle=\langle
\!\langle{f},{\varphi}\rangle\!\rangle\label{eq:ell-ds}%
\end{equation}
holds in $V^{\prime}$.
\end{definition}

As with evolution equations, we fix an orthonormal basis $\mathfrak{U}$ in
${\mathcal{U}}$ and use (\ref{eq:IWPG}) to rewrite (\ref{eq:ell}) as
\begin{equation}
{\mathbf{A}}u+\left(  {\mathbf{M}}u\right)  \diamond\dot{W}=f,
\label{eq:ell-b}%
\end{equation}
where
\begin{equation}
\mathbf{M}u\diamond\dot{W}:=\sum_{k\geq1} \mathbf{M}_{k}u\diamond\xi_{k}.
\label{eq:operM}%
\end{equation}
Taking $\varphi=\xi_{\alpha}$ in (\ref{eq:ell-ds}) and using relation
(\ref{eq:IBasG}) we conclude, as in Theorem \ref{th: equivalence1}, that
$u=\sum_{\alpha\in{\mathcal{J}}}u_{\alpha}\xi_{\alpha}$ is a solution of
equation (\ref{eq:ell}) if and only if $u_{\alpha}$ satisfies
\begin{equation}
{\mathbf{A}}u_{\alpha}+\sum_{k\geq1}\sqrt{\alpha_{k}}\;{\mathbf{M}}%
_{k}u_{\alpha-\varepsilon_{k}}=f_{\alpha} \label{eq:ell-S}%
\end{equation}
in the normal triple $(V,H,V^{\prime})$. This system of equation is
lower-triangular and can be solved by induction on $|\alpha|$.

The following example illucidates the limitations on the \textquotedblleft
quality\textquotedblright\ of the solution of equation (\ref{eq:ell}).

\begin{example}
\label{ex2} Consider equation
\begin{equation}
\label{eq:e2}u=1+u\diamond\xi.
\end{equation}
Similar to Example \ref{ex1}, we write $u=\sum_{n\geq0}u_{(n)}H_{n}(\xi
)/\sqrt{n!}$, where $H_{n}$ is Hermite polynomial of order $n$
(\ref{eq:hermite}). Then (\ref{eq:ell-S}) implies $u_{(n)}=I_{(n=0)}+\sqrt
{n}u_{(n-1)}$ or $u_{(0)}=1$, $u_{(n)}=\sqrt{n!}$, $n\geq1$, or $u=1+\sum
_{n\geq1}H_{n}(\xi)$. Clearly, the series does not converge in $L_{2}%
({\mathbb{F}})$, but does converge in $({\mathcal{S}})_{-1,q}$ for every $q<0$
(see (\ref{eq:S-weight})). As a result, even a simple stationary equation
(\ref{eq:e2}) can be solved only in weighted spaces.
\end{example}

\begin{theorem}
\label{th:ell} Consider equation (\ref{eq:ell-b}) in which $f\in
\bar{{\mathcal{R}}} L_{2}({\mathbb{F}};V^{\prime})$ for some $\bar
{{\mathcal{R}}}$.

Assume that the deterministic equation ${\mathbf{A}} U =F$ is uniquely
solvable in the normal triple $(V,H,V^{\prime})$, that is, for every $F\in
V^{\prime}$, there exists a unique solution $U={\mathbf{A}}^{-1}F\in V$ so
that $\|U\|_{V}\leq C_{A}\|F\|_{V^{\prime}}$. Assume also that each
${\mathbf{M}}_{k}$ is a bounded linear operator from $V$ to $V^{\prime}$ so
that, for all $v\in V$
\begin{equation}
\label{eq:op-norm-ell}\|{\mathbf{A}}^{-1}{\mathbf{M}}_{k} v\|_{V}\leq
C_{k}\|v\|_{V},
\end{equation}
with $C_{k}$ independent of $v$.

Then there exists an operator ${\mathcal{R}}$ and a unique solution
$u\in{\mathcal{R}} L_{2}({\mathbb{F}}; V)$ of (\ref{eq:evol-b}).
\end{theorem}

\textbf{Proof.} The argument is identical to the proof of Theorem
\ref{th:evol}. \endproof

\begin{remark}
\label{rm:ell} The assumption of the theorem about solvability of the
deterministic equation holds if the operator ${\mathbf{A}}$ satisfies
$\langle{\mathbf{A}} v,v\rangle\geq\kappa\|v\|^{2}_{V}$ for every $v\in V,$
with $\kappa>0$ independent of $v$.
\end{remark}

An analog of Theorem \ref{th:evol-nr} exists if $f$ is non-random. With no
time variable, we introduce the following notation to write multiple integrals
in the time-independent setting:
\[
{\delta}_{{\mathbf{B}}}^{(0)}(\eta)=\eta,\ {\delta}_{{\mathbf{B}}}^{(n)}%
(\eta)= {\delta}({\mathbf{B}}{\delta}_{{\mathbf{B}}}^{(n-1)}(\eta)), \ \eta
\in{\mathcal{R}} L_{2}({\mathbb{F}};V),
\]
where $\mathbf{B}$ is a bounded linear operator from $V$ to $V\otimes
\mathcal{U}$.

\begin{theorem}
\label{th:stat-nr} Under the assumptions of Theorem \ref{th:ell}, if $f$ is
non-random, then the following holds:

\begin{enumerate}
\item the coefficient $u_{\alpha}$, corresponding to the multi-index $\alpha$
with $|\alpha|=n\geq1$ and the characteristic set $K_{\alpha}= \{k_{1},
\ldots, k_{n}\}$, is given by
\begin{equation}
\label{coef-stat}u_{\alpha}=\frac{1}{\sqrt{\alpha!}}\sum_{\sigma
\in{\mathcal{P}}_{n}} {\mathbf{B}}_{k_{\sigma(n)}}\cdots{\mathbf{B}%
}_{k_{\sigma(1)}}u_{(0)},
\end{equation}
where

\begin{itemize}
\item ${\mathcal{P}}_{n}$ is the permutation group of the set $(1,\ldots, n)$;

\item $\mathbf{B}_{k}=-{\mathbf{A}}^{-1}{\mathbf{M}}_{k}$;

\item $u_{(0)}={\mathbf{A}}^{-1} f$.
\end{itemize}

\item the operator ${\mathcal{R}}$ can be defined by the weights $r_{\alpha}$
in the form
\begin{equation}
\label{eq:r-alpha-st}r_{\alpha}=\frac{q^{\alpha}}{\sqrt{|\alpha|!}},\ \mathrm{where\ } q^{\alpha}=\prod_{k=1}^{\infty}q_{k}^{\alpha_{k}},
\end{equation}
where the numbers $q_{k}, \, k\geq1$ are chosen so that $\sum_{k\geq1}
q_{k}^{2}C^{2}_{k} <1$, and $C_{k}$ are defined in (\ref{eq:op-norm-ell}).

\item With $r_{\alpha}$ and $q_{k}$ defined by (\ref{eq:r-alpha-st}),
\begin{equation}
\label{multiint-st1}\sum_{|\alpha|=n} q^{\alpha} u_{\alpha}\xi_{\alpha}=
{\delta}_{{\overline{{\mathbf{B}}}}}^{(n)}({\mathbf{A}}^{-1}f),
\end{equation}
where ${\overline{{\mathbf{B}}}}=-(q_{1}{\mathbf{A}}^{-1}{\mathbf{M}}_{1},
q_{2}{\mathbf{A}}^{-1}{\mathbf{M}}_{2},\ldots)$, and
\begin{equation}
\label{multiint-st2}{\mathcal{R}} u={\mathbf{A}}^{-1}f+\sum_{n\geq1}\frac{1}{\sqrt{n!}} \: {\delta}_{{\overline{{\mathbf{B}}}}}^{(n)}({\mathbf{A}%
}^{-1}f),
\end{equation}

\end{enumerate}
\end{theorem}

\textbf{Proof}. While the proofs of Theorems \ref{th:evol-nr} and
\ref{th:stat-nr} are similar, the complete absence of time makes equation
(\ref{eq:ell-b}) different from either (\ref{eq:evol-b}) or anything
considered in \cite{LR_AP}. Accordingly, we present a complete proof.

Define ${\widetilde{u}}_{\alpha}=\sqrt{\alpha!}\,u_{\alpha}$. If $f$ is
deterministic, then ${\widetilde{u}}_{(0)}={\mathbf{A}}^{-1}f$ and, for
$|\alpha|\geq1$,
\[
{\mathbf{A}}{\widetilde{u}}_{\alpha}+\sum_{k\geq1} \alpha_{k}{\mathbf{M}}_{k}
{\widetilde{u}}_{\alpha-\varepsilon_{k}}=0,
\]
or
\[
{\widetilde{u}}_{\alpha}=\sum_{k\geq1} \alpha_{k}{\mathbf{B}}_{k}%
{\widetilde{u}}_{\alpha-\varepsilon_{k}}= \sum_{k\in K_{\alpha}}{\mathbf{B}%
}_{k}{\widetilde{u}}_{\alpha-\varepsilon_{k}},
\]
where $K_{\alpha}=\{k_{1},\ldots, k_{n}\}$ is the characteristic set of
$\alpha$ and $n=|\alpha|$. By induction on $n$,
\[
{\widetilde{u}}_{\alpha}=\sum_{\sigma\in{\mathcal{P}}_{n}} {\mathbf{B}%
}_{k_{\sigma(n)}}\cdots{\mathbf{B}}_{k_{\sigma(1)}}u_{(0)},
\]
and (\ref{coef-stat}) follows.

Next, define
\[
U_{n}=\sum_{|\alpha|=n}q^{\alpha}u_{\alpha}\xi_{\alpha}, \ n\geq0.
\]
Let us first show that, for each $n\geq1$, $U_{n}\in L_{2}({\mathbb{F}}; V)$.
By (\ref{coef-stat}) we have
\begin{equation}
\label{eq:ell-sp-m}\|u_{\alpha}\|_{V}^{2}\leq C_{A}^{2} \frac{(|\alpha|!)^{2}%
}{\alpha!}\|f\|_{V^{\prime}}^{2}\, \prod_{k\geq1}C_{k}^{\alpha_{k}}.
\end{equation}
By (\ref{multinomial}),
\begin{equation*}
\begin{split}
\sum_{|\alpha|=n}q^{2\alpha}\|u_{\alpha}\|_{V}^{2} &
\leq C_{A}^{2}
\|f\|_{V^{\prime}}^{2}\,n!\,\sum_{|\alpha|=n}\left( \frac{n!}{\alpha!} \prod_{k\geq1} (C_{k}q_{k})^{2\alpha_{k}}\right)\\
&=C_{A}^{2}\|f\|_{V^{\prime}}^{2}\,n!\, \left(  \sum_{k\geq1}C^{2}_{k}q^{2}_{k}\right)^{n} < \infty,
\end{split}
\end{equation*}
because of the selection of $q_{k}$, and so $U_{n}\in L_{2}({\mathbb{F}}; V)$.
If the weights $r_{\alpha}$ are defined by (\ref{eq:r-alpha-st}), then
\[
\sum_{\alpha\in{\mathcal{J}}} r_{\alpha}^{2} \|u\|_{V}^{2}= \sum_{n\geq0}%
\sum_{|\alpha|=n}r_{\alpha}^{2} \|u\|_{V}^{2} \leq C_{A}^{2}\|f\|_{V^{\prime}}^{2}\sum_{n\geq 0}\left(  \sum_{k\geq1}C^{2}_{k}q^{2}_{k}\right)  ^{n} < \infty,
\]
because of the assumption $\sum_{k\geq1}C^{2}_{k}q^{2}_{k}< 1$.

Since (\ref{multiint-st2}) follows directly from (\ref{multiint-st1}), it
remains to establish (\ref{multiint-st1}), that is,
\begin{equation}
\label{st-aux11}U_{n}={\delta}_{{\overline{{\mathbf{B}}}}}(U_{n-1}),\ n\geq1.
\end{equation}
For $n=1$ we have
\[
U_{1}=\sum_{k\geq1} q_{k}u_{\varepsilon_{k}}\xi_{k}= \sum_{k\geq1}
{\overline{{\mathbf{B}}}}_{k}u_{(0)}\xi_{k}={\delta}_{{\overline{{\mathbf{B}}%
}}}(U_{0}),
\]
where the last equality follows from (\ref{eq:IWPG}). More generally, for
$n>1$ we have by definition of $U_{n}$ that
\[
(U_{n})_{\alpha}=
\begin{cases}
q^{\alpha}u_{\alpha}, & \mathrm{if} \ |\alpha|=n,\\
0, & \mathrm{otherwise}.
\end{cases}
\]
From the equation
\[
q^{\alpha}{\mathbf{A}} u_{\alpha}+ \sum_{k\geq1}q_{k}\sqrt{\alpha_{k}}\,
{\mathbf{M}}_{k} q^{\alpha-\varepsilon_{k}} u_{\alpha-\varepsilon_{k}}=0
\]
we find
\[%
\begin{split}
(U_{n})_{\alpha}  &  =
\begin{cases}
\displaystyle \sum_{k\geq1}\sqrt{\alpha_{k}}\, q_{k} {\mathbf{B}}_{k}
q^{\alpha-\varepsilon_{k}} u_{\alpha-\varepsilon_{k}}, & \mathrm{if}
\ |\alpha|=n,\\
0, & \mathrm{otherwise}.
\end{cases}
\\
&  =\sum_{k\geq1}\sqrt{\alpha_{k}}\ {\overline{{\mathbf{B}}}}_{k}%
(U_{n-1})_{\alpha-\varepsilon_{k}},
\end{split}
\]
and then (\ref{st-aux11}) follows from (\ref{eq:IBasG}). Theorem
\ref{th:stat-nr} is proved. \endproof

Here is another result about solvability of (\ref{eq:ell-b}), this time with
random $f$. We use the space $({\mathcal{S}})_{\rho, q}$, defined by the
weights (\ref{eq:S-weight}).

\begin{theorem}
\label{th:ell-sp2} In addition to the assumptions of Theorem \ref{th:ell}, let
$C_{A}\leq1$ and $C_{k}\leq1$ for all $k$. If $f\in({\mathcal{S}})_{-1,-\ell
}(V^{\prime})$ for some $\ell>1$, then there exists a unique solution
$u\in({\mathcal{S}})_{-1,-\ell-4}(V)$ of (\ref{eq:ell-b}) and
\begin{equation}
\label{eq:ell-sp3}\|u\|_{({\mathcal{S}})_{-1,-\ell-4}(V)} \leq C(\ell
)\|f\|_{({\mathcal{S}})_{-1,-\ell}(V^{\prime})}.
\end{equation}

\end{theorem}

\textbf{Proof.} Denote by $u(g;\gamma)$, $\gamma\in{\mathcal{J}}$, $g\in
V^{\prime}$, the solution of (\ref{eq:ell-b}) with $f_{\alpha}=gI_{(\alpha
=\gamma)}$, and define $\bar{u}_{\alpha}=(\alpha!)^{-1/2}u_{\alpha}$. Clearly,
$u_{\alpha}(g,\gamma)=0$ if $|\alpha|<|\gamma|$ and so
\begin{equation}
\label{eq:ell-sp4}\sum_{\alpha\in{\mathcal{J}}} \|u_{\alpha}(f_{\gamma}%
;\gamma)\|_{V}^{2}r_{\alpha}^{2}= \sum_{\alpha\in{\mathcal{J}}}\|u_{\alpha
+\gamma}(f_{\gamma};\gamma)\|_{V}^{2}r_{\alpha+\gamma}^{2}.
\end{equation}
It follows from (\ref{eq:ell-S}) that
\begin{equation}
\label{eq:ell-sp5}\bar{u}_{\alpha+\gamma}(f_{\gamma};\gamma)=\bar{u}_{\alpha
}\big(f_{\gamma}(\gamma!)^{-1/2}; (0)\big).
\end{equation}
Now we use (\ref{eq:ell-sp-m}) to conclude that
\begin{equation}
\label{eq:ell-sp6}\|\bar{u}_{\alpha+\gamma}(f_{\gamma};\gamma)\|_{V}\leq
\frac{|\alpha|!}{\sqrt{\alpha!\gamma!}} \|f\|_{V^{\prime}}.
\end{equation}
Coming back to (\ref{eq:ell-sp4}) with $r_{\alpha}^{2}=(\alpha!)^{-1}%
(2{\mathbb{N}})^{(-\ell-4)\alpha}$ and using inequality (\ref{HOUZ-ineq}) we
find:
\[
\|u(f_{\gamma};\gamma)\|_{({\mathcal{S}})_{-1,-\ell-4}(V)}\leq C(\ell
)(2{\mathbb{N}})^{-2\gamma}\, \frac{\|f_{\gamma}\|_{V^{\prime}}}
{(2{\mathbb{N}})^{(\ell/2) \gamma}\sqrt{\gamma!}},
\]
where
\[
C(\ell)=\left(  \sum_{\alpha\in{\mathcal{J}}}\left(  \frac{|\alpha|!}{\alpha
!}\right)  ^{2} (2\mathbb{N})^{(-\ell-4)\alpha}\right)  ^{1/2};
\]
(\ref{eq:HOUZ1}) and (\ref{HOUZ-ineq}) imply $C(\ell)<\infty$. Then
(\ref{eq:ell-sp3}) follows by the triangle inequality after summing over all
$\gamma$ and using the Cauchy-Schwartz inequality. \endproof

\begin{remark}
\label{rm:ell-sp1} Example \ref{ex2}, in which $f\in({\mathcal{S}})_{0,0}$ and
$u\in({\mathcal{S}})_{-1,q}$, $q<0$, shows that, while the results of Theorem
\ref{th:ell-sp2} are not sharp, a bound of the type $\Vert u\Vert
_{({\mathcal{S}})_{\rho,q}(V)}\leq C\Vert f\Vert_{({\mathcal{S}})_{\rho,\ell
}(V^{\prime})}$ is, in general, impossible if $\rho>-1$ or $q\geq\ell$.
\end{remark}

\subsection{ Convergence to Stationary Solution}

Let $(V,H,V^{\prime})$ be a normal triple of Hilbert spaces. Consider
equation
\begin{equation}
\dot{u}(t)=({\mathbf{A}}u(t)+f(t))+{\mathbf{M}}_{k}u(t)\diamond\xi_{k},
\label{eq:conv1}%
\end{equation}
where the operators ${\mathbf{A}}$ and ${\mathbf{M}}_{k}$ do not depend on
time, and assume that there exists an $f^{\ast}\in{\mathcal{R}}L_{2}%
(\mathbb{F};H)$ such that $\lim_{t\rightarrow\infty}\Vert f(t)-f^{\ast}%
\Vert_{{\mathcal{R}}L_{2}({\mathbb{F}};H)}=0$. The objective of this section
is to study convergence, as $t\rightarrow+\infty$, of the solution of
(\ref{eq:conv1}) to the solution $u^{\ast}$ of the stationary equation
\begin{equation}
-{\mathbf{A}}u^{\ast}=f^{\ast}+{\mathbf{M}}_{k}u^{\ast}\diamond\xi_{k}.
\label{eq:conv2}%
\end{equation}

\begin{theorem}
\label{th:conv} Assume that

\begin{enumerate}
\item[(\textbf{C1})] Each ${\mathbf{M}}_{k}$ is a bounded linear operator from
$H$ to $H$, and ${\mathbf{A}}$ is a bounded linear operator from $V$ to
$V^{\prime}$ with the property
\begin{equation}
\label{eq:conv3}\langle{\mathbf{A}} v,v\rangle+ \kappa\|v\|^{2}_{V}%
\leq-c\|v\|_{H}^{2}%
\end{equation}
for every $v\in V$, with $\kappa>0$ and $c>0$ both independent of $v$.

\item[(\textbf{C2})] $f\in\bar{{\mathcal{R}}} L_{2}({\mathbb{F}};{\mathcal{H}%
})$ and there exists an $f^{*}\in\bar{{\mathcal{R}}} L_{2}({\mathbb{F}};H)$
such that \newline$\lim_{t\to+\infty}\|f(t)-f^{*}\|_{\bar{{\mathcal{R}}}
L_{2}({\mathbb{F}};H)}$.
\end{enumerate}

Then, for every $u_{0}\in\bar{{\mathcal{R}}} L_{2}({\mathbb{F}};H),$ there
exists an operator ${{\mathcal{R}}}$ so that

\begin{enumerate}
\item There exists a unique solution $u\in{{\mathcal{R}}} L_{2}({\mathbb{F}%
};{\mathcal{V}})$ of (\ref{eq:conv1}),

\item There exists a unique solution $u^{*}\in{{\mathcal{R}}} L_{2}%
({\mathbb{F}};V)$ of (\ref{eq:conv2}), and

\item The following convergence holds:
\begin{equation}
\lim_{t\rightarrow+\infty}\Vert u(t)-u^{\ast}\Vert_{{{\mathcal{R}}}%
L_{2}({\mathbb{F}};H)}=0. \label{eq:conv4}%
\end{equation}

\end{enumerate}
\end{theorem}

\textbf{Proof} (1) Existence and uniqueness of the solution of (\ref{eq:conv1}%
) follow from Theorem \ref{th:evol} and Remark \ref{rm:evol}.

(2) Existence and uniqueness of the solution of (\ref{eq:conv2}) follow from
Theorem \ref{th:ell} and Remark \ref{rm:ell}.

(3) The proof of (\ref{eq:conv4}) is based on the following result.

\begin{lemma}
\label{lm:conv} Assume that the operator ${\mathbf{A}}$ satisfies
(\ref{eq:conv3}) and $F=F(t)$ is a deterministic function such that
$\lim_{t\to+\infty}\|F(t)\|_{H}=0$. Then, for every $U_{0}\in H$, the solution
$U=U(t)$ of the equation $U(t)=U_{0}+\int_{0}^{t}{\mathbf{A}} U(s)ds +
\int_{0}^{t} F(s)ds$ satisfies $\lim_{t\to+\infty}\|U(t)\|_{H}=0$.
\end{lemma}

\textbf{Proof.} If $\Phi=\Phi_{t}$ is the semi-group generated by the operator
${\mathbf{A}}$ (which exists because of (\ref{eq:conv3})), then
\[
U(t)=\Phi_{t}U_{0}+\int_{0}^{t}\Phi_{t-s}F(s)ds.
\]
Condition (\ref{eq:conv3}) implies $\|\Phi_{t} U_{0}\|_{H} \leq e^{-ct}%
\|U_{0}\|_{H}$, and then
\[
\|U(t)\|_{H}\leq e^{-ct}\|U_{0}\|_{H}+\int_{0}^{t} e^{-c(t-s)}\|F(s)\|_{H}ds.
\]
The convergence of $\|U(t)\|_{H}$ to zero now follows from the Toeplitz lemma
(see Lemma A.2 in Appendix). Lemma \ref{lm:conv} is proved. \endproof

To complete the proof of Theorem \ref{th:conv}, we define $v_{\alpha
}(t)=u_{\alpha}(t)-u_{\alpha}^{\ast}$ and note that
\[
\dot{v}_{\alpha}(t)={\mathbf{A}}v_{\alpha}(t)+(f_{\alpha}(t)-f_{\alpha}^{\ast
})+\sum_{k}\sqrt{\alpha_{k}}\, {\mathbf{M}}_{k}v_{\alpha-\varepsilon_{k}}.
\]
By Theorem \ref{th:ell}, $u_{\alpha}^{\ast}\in V$ and so $v_{\alpha}(0)\in H$
for every $\alpha\in{\mathcal{J}}$. By Lemma \ref{lm:conv}, $\lim
_{t\rightarrow+\infty}\Vert v_{(0)}(t)\Vert_{H}=0$. Using induction on
$|\alpha|$ and the inequality $\Vert{\mathbf{M}}_{k}v_{\alpha-\varepsilon_{k}%
}(t)\Vert_{H}\leq c_{k}\Vert v_{\alpha-\varepsilon_{k}}(t)\Vert_{H}$, we
conclude that $\lim_{t\rightarrow+\infty}\Vert v_{\alpha}(t)\Vert_{H}=0$ for
every $\alpha\in{\mathcal{J}}$. Since $v_{\alpha}\in{\mathbf{C}}((0,T);H)$ for
every $T$, it follows that $\sup_{t\geq0}\Vert v_{\alpha}(t)\Vert_{H}<\infty$.
Define the operator ${\mathcal{R}}$ on $L_{2}({\mathbb{F}})$ so that
${\mathcal{R}}\xi_{\alpha}=r_{\alpha}\xi_{\alpha}$, where
\[
r_{\alpha}=\frac{(2{\mathbb{N}})^{-\alpha}}{1+\sup\limits_{t\geq0}\Vert
v_{\alpha}(t)\Vert_{H}}.
\]
Then (\ref{eq:conv4}) follows by the dominated convergence theorem.

Theorem \ref{th:conv} is proved. \endproof

\section{Bilinear parabolic and elliptic SPDEs}

\label{sec5} \setcounter{equation}{0} \setcounter{theorem}{0}

Let $G$ be a smooth bounded domain in ${\mathbb{R}}^{d}$ and $\{h_{k}%
,\;k\geq1\}$, an orthonormal basis in $L_{2}(G)$. We assume that
\begin{equation}
\sup_{x\in G}|h_{k}(x)|\leq c_{k},\ k\geq1. \label{eq:wn-base}%
\end{equation}
A space white noise on $L_{2}(G)$ is a formal series
\begin{equation}
{\dot{W}}(x)=\sum_{k\geq1}h_{k}(x)\xi_{k}, \label{eq:wn-wn}%
\end{equation}
where $\xi_{k},\ k\geq1,$ are independent standard Gaussian random variables.

\subsection{Dirichlet Problem for parabolic SPDE of the Second Order}

\ Consider the following equation:
\begin{equation}%
\begin{split}
u_{t}(t,x)  &  =a_{ij}(x)D_{i}D_{j}u(t,x)+b_{i}(x)D_{i}%
u(t,x)+c(x)u(t,x)+f(t,x)\\
&  +(\sigma_{i}(x)D_{i}u(t,x)+\nu(x)u(t,x)+g(t,x))\diamond{\dot{W}%
}(x),\ 0<t\leq T,\ x\in G,
\end{split}
\label{eq:wn-evol}%
\end{equation}
with zero boundary conditions and some initial condition $u(0,x)=u_{0}(x)$;
the functions $a_{ij},\,b_{i},\,c,\,f,\,\sigma_{i},\,\nu,\,g,$ and $u_{0}$ are
non-random. In (\ref{eq:wn-evol}) and in similar expressions below we assume
summation over the repeated indices. Let $(V,H,V^{\prime})$ be the normal
triple with ${V=\overset{\circ}{H}{}_{2}^{1}(G)}$, $H=L_{2}(G)$, $V^{\prime
}=H_{2}^{-1}(G)$. In view of (\ref{eq:wn-wn}), equation (\ref{eq:wn-evol}) is
a particular case of equation (\ref{eq:evol-b}) so that
\begin{equation}
{\mathbf{A}}u=a_{ij}(x)D_{i}D_{j}u+b_{i}(x)D_{i}u+c(x)u,\ {\mathbf{M}}%
_{k}u=(\sigma_{i}(x)D_{i}u+\nu(x)u)h_{k}(x), \label{eq:wn-oper}%
\end{equation}
and $f(t,x)+g(t,x)\diamond{\dot{W}}(x)$ is the free term.

We make the following assumptions about the coefficients:

\begin{enumerate}
\item[\textbf{D1}] The functions $a_{ij}$ are Lipschitz continuous in the
closure $\bar{G}$ of $G$, and the functions $b_{i},\;c,\;\sigma_{i},\;\nu$ are
bounded and measurable in $\bar{G}$.

\item[\textbf{D2}] There exist positive numbers $A_{1},A_{2}$ so that
$A_{1}|y|^{2}\leq a_{ij}(x)y_{i}y_{j}\leq A_{2}|y|^{2}$ for all $x\in\bar{G}$
and $y\in{\mathbb{R}}^{d}$.
\end{enumerate}

Given a $T>0$, recall the notation ${\mathcal{V}}=L_{2}((0,T);V)$ and
similarly for ${\mathcal{H}}$ and ${\mathcal{V}}^{\prime}$ (see
(\ref{sp-notation})).

\begin{theorem}
\label{th:wn-evol} Under the assumptions \textbf{D1} and \textbf{D2}, if
$f\in{\mathcal{V}}^{\prime}$, $g\in{\mathcal{H}}$, $u_{0}\in H$, then there
exists an $\ell>1$ and a number $C>0$, both independent of $u_{0},f,g$, so
that $u\in{\mathcal{R}}L_{2}({\mathbb{F}};{\mathcal{V}})$ and
\begin{equation}
\Vert u\Vert_{{\mathcal{R}}L_{2}({\mathbb{F}};{\mathcal{V}})}\leq
C\cdot\big(\Vert u_{0}\Vert_{H}+\Vert f\Vert_{{\mathcal{V}}^{\prime}}+\Vert
g\Vert_{{\mathcal{H}}}\big), \label{eq:wn-evol-th}%
\end{equation}
where the operator ${\mathcal{R}}$ is defined by the weights
\begin{equation}
r_{\alpha}^{2}=c^{-2\alpha}(|\alpha|!)^{-1}(2{\mathbb{N}})^{-2\ell\alpha}
\label{eq:wn-weight}%
\end{equation}
and $c^{\alpha}=\prod_{k}c_{k}^{\alpha_{k}}$, with $c_{k}$ from
(\ref{eq:wn-base}); the number $\ell$ in general depends on $T$.
\end{theorem}

\textbf{Proof.} We derive the result from Theorem \ref{th:evol-nr}. Consider
the deterministic equation $\dot{U}(t)={\mathbf{A}} U(t)+F$. Assumptions
\textbf{D1} and \textbf{D2} imply that there exists a unique solution of this
equation in the normal triple $(V,H,V^{\prime})$, and the solution satisfies
\begin{equation}
\label{eq:wn-det}\sup_{0<t<T}\|U(t)\|_{H}+\|U\|_{{\mathcal{V}}}\leq
C\cdot\big(\|U(0)\|_{H}+\|F\|_{{\mathcal{V}}^{\prime}}\big),
\end{equation}
where the number $C$ depends on $T$ and the operator ${\mathbf{A}}$. Moreover,
(\ref{eq:wn-base}) implies that (\ref{Ck-def}) holds with $C_{k}=C_{0}c_{k}$
for some positive number $C_{0}$ independent of $k$, but possibly depending on
$T$.

To proceed, let us assume first that $g=0$. Then the statement of the theorem
follows directly from Theorem \ref{th:evol-nr} if we take in
(\ref{r-alpha-evol}) $q_{k}=c_{k}{}^{-1}\left(  2k\right)  ^{-\ell}$ with
sufficiently large $\ell$.

It now remains to consider the case $g\not =0$ and $f=u_{0}=0$. Even though
$g$ is non-random, $g\xi_{k}$ is, and therefore a direct application of
Theorem \ref{th:evol-nr} is not possible. Instead, let us look more closely at
the corresponding equations for $u_{\alpha}$. For $\alpha=(0)$,
\[
u_{(0)}(t)=\int_{0}^{t} {\mathbf{A}} u_{(0)}(s)ds,
\]
which implies $u_{(0)}(t)=0$ for all $t$. For $\alpha=\varepsilon_{k}$,
\[
{u}_{\varepsilon_{k}}(t)=\int_{0}^{t} {\mathbf{A}} {u}_{\varepsilon_{k}}(s)ds+
h_{k}\int_{0}^{t}g(s)ds,
\]
or
\[
u_{\varepsilon_{k}}(t)=\int_{0}^{t}\Phi_{t-s}h_{k}g(s)ds,
\]
so that
\begin{equation}
\label{eq:wn-aux00}\|u_{\varepsilon_{k}}\|_{{\mathcal{V}}}\leq C_{0}%
c_{k}\|g\|_{{\mathcal{H}}}.
\end{equation}
If $|\alpha|>1$, then
\[
u_{\alpha}(t)=\int_{0}^{t}{\mathbf{A}} u_{\alpha}(s)ds+\sum_{k\geq1}
\sqrt{\alpha_{k}} \,{\mathbf{M}}_{k} u_{\alpha-\varepsilon_{k}},
\]
which is the same as (\ref{eq:evol-S-det1}). In particular, if $|\alpha|=2$
and $\{i,j\}$ is the characteristic set of $\alpha$, then
\[
u_{\alpha}(t)=\frac{1}{\sqrt{\alpha!}} \int_{0}^{t} \Phi_{t-s}\left(
{\mathbf{M}}_{i} u_{\varepsilon_{j}}(s)+{\mathbf{M}}_{j} u_{\varepsilon_{i}%
}(s) \right)  ds.
\]
More generally, by analogy with (\ref{aux-evol-nr000}), if $|\alpha|=n>2$ and
$\{k_{1}, \ldots, k_{n}\}$ is the characteristic set of $\alpha$, then
\[
u_{\alpha}(t)=\frac{1}{\sqrt{\alpha!}}\sum_{\sigma\in{\mathcal{P}}_{n}}
\int_{0}^{t}\int_{0}^{s_{n}}\ldots\int_{0}^{s_{3}} \Phi_{t-s_{n}}{\mathbf{M}%
}_{k_{\sigma(n)}}\cdots\Phi_{s_{3}-s_{2}} {\mathbf{M}}_{k_{\sigma(2)}%
}u_{\varepsilon_{\sigma(1)}}(s_{2})ds_{2}\ldots ds_{n}.
\]
By the triangle inequality and (\ref{eq:wn-aux00}),
\[
\|u_{\alpha}\|_{{\mathcal{V}}}\leq\frac{|\alpha|! C_{0}^{|\alpha|}c^{\alpha}%
}{\sqrt{\alpha!}}\, \|g\|_{{\mathcal{H}}},
\]
and then (\ref{eq:wn-evol-th}) follows from \eqref{HOUZ-ineq} if $\ell$ is
sufficiently large.

This completes the proof of Theorem \ref{th:wn-evol}. \endproof

\begin{theorem}
\label{th:wn-conv} In addition to \textbf{D1} and \textbf{D2}, assume that

\begin{enumerate}
\item $\sigma_{i}=0$ for all $i$;

\item the operator ${\mathbf{A}}$ in $G$ with zero boundary conditions
satisfies (\ref{eq:conv3}).
\end{enumerate}

If there exist functions $f^{*}$ and $g^{*}$ from $H$ so that
\begin{equation}
\label{eq:wn-conv1}\lim_{t\to+\infty}\left(  \|f(t)-f^{*}\|_{H}+
\|g(t)-g^{*}\|_{H}\right)  =0,
\end{equation}
then the solution $u$ of equation (\ref{eq:wn-evol}) satisfies
\begin{equation}
\label{eq:wn-conv2}\lim_{t\to+\infty}\|u(t)-u^{*}\|_{{\mathcal{R}}
L_{2}({\mathbb{F}};H)}=0,
\end{equation}
where the operator ${\mathcal{R}}$ is defined by the weights
(\ref{eq:wn-weight}) and $u^{*}$ is the solution of the stationary equation
\begin{equation}
\label{eq:wn-ste}%
\begin{split}
a_{ij}(x)D_{i}D_{j}u^{*}(x)  &  +b_{i}(x)D_{i}u^{*}(x)+c(x)u^{*}(x)+f^{*}(x)\\
&  \!\!+(\nu(x)u^{*}(x)+g^{*}(x))\diamond{\dot{W}}(x)=0, \ x\in
G;\ u|_{\partial G}=0.
\end{split}
\end{equation}

\end{theorem}

\textbf{Proof.} This follows from Theorem \ref{th:conv}. \endproof

\begin{remark}
The operator ${\mathbf{A}}$ satisfies (\ref{eq:conv3}) if, for example, each
$a_{ij}$ is twice continuously differentiable in $\bar{G}$, each $b_{i}$
continuously differentiable in $\bar{G}$, and
\begin{equation}
\inf_{x\in\bar{G}}c(x)-\sup_{x\in\bar{G}}(|D_{i}D_{j}a_{ij}(x)|+|D_{i}%
b_{i}(x)|)\geq\varepsilon>0; \label{eq:C3}%
\end{equation}
this is verified directly using integration by parts.
\end{remark}

\subsection{Elliptic SPDEs of the full second order}

Consider the following Dirichlet problem:%
\begin{equation}%
\begin{array}
[c]{c}%
- D_{i}\Big(a_{ij}\left(  x\right)  D_{j} u\left(  x\right)  \Big) +\\
D_{i}\Big(\sigma_{ij}\left(  x\right)  D_{j}\left(  u\left(  x\right)
\right)  \Big) \diamond\dot{W}\left(  x\right)  =f\left(  x\right)  ,\text{
}x\in G,\\
u_{|\partial G}=0,
\end{array}
\label{eq:wn-ell}%
\end{equation}
where $\dot{W}$ is the space white noise (\ref{eq:wn-wn}). Assume that the
functions $a_{ij},\,\sigma_{ij},\,f,\,$ and $g$ are non-random. Recall that
according to our summation convention, in (\ref{eq:wn-ell}) and in similar
expressions below we assume summation over the repeated indices.

We make the following assumptions:

\begin{description}
\item[E1] The functions $a_{ij}=a_{ij}(x)$ and $\sigma_{ij}=\sigma_{ij}(x)$
are measurable and bounded in the closure $\bar{G}$ of $G$.

\item[E2] There exist positive numbers $A_{1},A_{2}$ so that $A_{1}|y|^{2}\leq
a_{ij}(x)y_{i}y_{j}\leq A_{2}|y|^{2}$ for all $x\in\bar{G}$ and $y\in
{\mathbb{R}}^{d}$.

\item[E3] The functions $h_{k}$ in (\ref{eq:wn-wn}) are bounded and Lipschitz continuous.
\end{description}

Clearly, equation (\ref{eq:wn-ell}) is a particular case of equation
(\ref{eq:ell-b}) with
\begin{equation}
{\mathbf{A}}u(x):=-D_{i}\Big( a_{ij}\left(  x\right)  D_{j}u\left(  x\right)
\Big) \ \label{eq:wn-drift}%
\end{equation}
and
\begin{equation}
\ {\mathbf{M}}_{k}u(x):=h_{k}(x)\, D_{i}\Big( \sigma_{ij}\left(  x\right)
D_{j}u\left(  x\right)  \Big). \label{eq:wn-dif}%
\end{equation}
Assumptions \textbf{E1} and \textbf{E3} imply that each $\mathbf{M}_{k}$ is a
bounded linear operator from ${\overset{\circ}{H_{2}{}^{1}}(G)}$ to
$H_{2}^{-1}(G).$ Moreover, it is a standard fact that under the assumptions
\textbf{E1} and \textbf{E2} the operator $\mathbf{A}$ is an isomorphism from
$V$ onto $V^{\prime}$ (see e.g. \cite{LM}). Therefore, for every $k$ there
exists a positive number $C_{k}$ such that
\begin{equation}
\label{eq:fso}\left\Vert \mathbf{A}^{-1}M_{k}v\right\Vert _{V}
\leq C_{k}\left\Vert
v\right\Vert _{V},\ v\in V.
\end{equation}

\begin{theorem}
\label{th:wn-fso} Under the assumptions \textbf{E1} and \textbf{E2}, if $f\in
H_{2}^{-1}(G)$, then there exists a unique solution of the Dirichlet problem
(\ref{eq:wn-ell}) $u\in{\mathcal{R}}L_{2}({\mathbb{F}}; {\overset{\circ}{H}%
{}_{2}^{1}}(G))$ such that
\begin{equation}
\Vert u\Vert_{{\mathcal{R}}L_{2}({\mathbb{F}}; {\overset{\circ}{H}{}_{2}^{1}%
}(G))}\leq C\cdot\Vert f\Vert_{H_{2}^{-1}(G)}. \label{eq:wn-sta-th}%
\end{equation}
The weights $r_{\alpha}$ can be taken in the form
\begin{equation}
r_{\alpha}=\frac{q^{\alpha}}{\sqrt{|\alpha|!}},\ \mathrm{where\ }%
q^{\alpha}=\prod_{k=1}^{\infty}q_{k}^{\alpha_{k}}, \label{eq:r-alpha-exsta}%
\end{equation}
and the numbers $q_{k},\,k\geq1$ are chosen so that
$\sum_{k\geq1}C_{k}^{2}q_{k}^{2}<1$, with $C_{k}$ from (\ref{eq:fso}).
\end{theorem}

\textbf{Proof.} This follows from Theorem \ref{th:stat-nr}. \endproof

\begin{remark}
With an appropriate change of the boundary conditions, and with extra
regularity of the basis functions $h_{k}$, the results of Theorem
\ref{th:wn-fso} can be extended to stochastic elliptic equations of order
$2m$. The corresponding operators are
\begin{equation}
{\mathbf{A}}u=\left(  -1\right)  ^{m} D_{i_{1}}\cdots D_{i_{m}}\Big( a_{i_{1}%
\ldots i_{m}j_{1}\ldots j_{m}}\left(  x\right)  D_{j_{1}}\cdots D_{j_{m}}
u\left(  x\right)  \Big) \ \label{eq: wn-2m-drift}%
\end{equation}
and
\begin{equation}
\ {\mathbf{M}}_{k}u=h_{k}(x)\, D_{i_{1}}\cdots D_{i_{m}}\Big( \sigma
_{i_{1}\ldots i_{m}j_{1}\ldots j_{m}}\left(  x\right)  D_{j_{1}}\cdots
D_{j_{m}} u\left(  x\right)  \Big). \label{eq: wn-2m-dif}%
\end{equation}
Since $G$ is a smooth bounded domain, regularity of $h_{k}$ is not a problem:
we can take $h_{k}$ as the eigenfunctions of the Dirichlet Laplacian in $G$.
\end{remark}


\providecommand{\bysame}{\leavevmode\hbox to3em{\hrulefill}\thinspace}
\providecommand{\MR}{\relax\ifhmode\unskip\space\fi MR }
\providecommand{\MRhref}[2]{  \href{http://www.ams.org/mathscinet-getitem?mr=#1}{#2}
} \providecommand{\href}[2]{#2}

\textbf{Appendix}.

\textsc{A factorial inequality.}

\textbf{Lemma A.1.} \emph{For every multi-index $\alpha\in \mathcal{J}$,}
\[
|\alpha|! \leq \alpha! (2\mathbb{N})^{2\alpha}.
\tag{A1}
\]

\textbf{Proof.} Recall that, for $\alpha=(\alpha_1,\ldots, \alpha_k)\in
\mathcal{J}$,
$$
|\alpha|=\sum_{\ell=1}^k \alpha_{\ell}, \
\alpha!=\prod_{\ell=1}^k\alpha_{\ell}!,\
\mathbb{N}^{\alpha}=\prod_{\ell=1}^k \ell^{\alpha_{\ell}}.
$$
It is therefore clear that, if $|\alpha|=n$, then it is enough to
establish (A1) for $\alpha$ with $\alpha_{k}=0$ for $k\geq n+1$, because
a shift of a multi-index  entry to the right increases the right-hand side of
(A1) but does not change the left-had side. For example, if
$\alpha=(1,3,2,0,\ldots)$ and $\beta=(1,3,0,2,0,\ldots)$, then
$|\alpha|=|\beta|,\  \alpha!=\beta!, $ but
$\mathbb{N}^{\alpha}<\mathbb{N}^{beta}$.
Then
$$
4^n\geq \left(1+\frac{1}{2^2}+\ldots+\frac{1}{n^2}\right)^n =
\sum_{\alpha_1+ \ldots+ \alpha_n=n}
\frac{|\alpha|!}{\alpha!}\, \frac{1}{\mathbb{N}^{2\alpha}},
$$
where the equality follows by the multinomial formula. Since all the
term in the sum are non-negative, we get (A2).

\endproof

 The proof shows that inequality (A2) can be improved
by observing that $\sum_{k\geq1 }k^{-2}=\pi^2/6<2$. One can also
consider  $\sum_{k\geq 1} k^{-q}$ for some $1<q<2$.

\textsc{A version of the Toeplitz lemma.}

\textbf{Lemma A.2.} \emph{Assume that $f=f(t)$ is an integrable function and
$\lim_{t\to+\infty} |f(t)|=0$. Then, for every $c>0$, $\lim_{t\to+\infty}%
\int_{0}^{t}e^{-c(t-s)}f(s)ds=0$. }

\textbf{Proof.} Given $\varepsilon>0$, choose $T$ so that $|f(t)|<\varepsilon$
for all $t>T$. Then $|\int_{0}^{t}e^{-c(t-s)}f(s)ds|\leq e^{-ct}\int_{0}%
^{T}e^{cs} |f(s)|ds+ \varepsilon\int_{T}^{t}e^{-c(t-s)}ds$. Passing to the
limit as $t\to+\infty$, we find $\lim_{t\to+\infty}|\int_{0}^{t}%
e^{-c(t-s)}f(s)ds|\leq\varepsilon/c,$ which completes the proof.

 \endproof

\end{document}